\pdfoutput=1

\documentclass[10pt]{amsart}
\usepackage{graphicx}

\newtheorem{theorem}{Theorem}[section]

\theoremstyle{definition}

\begin{document}

\title[The Kuratowski covering conjecture]{The Kuratowski covering conjecture for graphs of order $<$ 10 for the nonorientable surfaces of genus 3 and 4}


\begin{abstract}

\noindent Kuratowski proved that a finite graph embeds in the plane 
if it does not contain a subdivision of either $K_5$ or $K_{3,3}$, 
called Kuratowski subgraphs.
A conjectured generalization of this result to all nonorientable surfaces 
says that a finite minimal forbidden subgraph for the nonorientable surface 
of genus $\tilde{g}$ can be written as the union of $\tilde{g}+1$ 
Kuratowski subgraphs such that 
the union of each pair of these fails to embed in the projective plane, 
the union of each triple of these fails to embed in the Klein bottle 
if $\tilde{g} \ge 2$, 
and the union of each triple of these fails to embed in the torus 
if $\tilde{g} \ge 3$.
We show that this conjecture is true for all minimal forbidden subgraphs 
of order $<$ 10 for the nonorientable surfaces of genus 3 and 4.

\end{abstract}

\author{Suhkjin Hur}

\address{
Department of Mathematics\\ 
The Ohio State University\\
231 West 18th Avenue \\
Columbus, OH 43210
}
\email{sjin@math.ohio-state.edu}





\maketitle

\section{Introduction}
\label{sec:intro}

We use the same terminology as in \cite{MT} unless otherwise specified.
$G_1 \vee G_2$ denotes the graph obtained by identifying 
one vertex of $G_1$ and one vertex of $G_2$.

Kuratowski \cite{K} showed that minimal forbidden subgraphs for the plane 
are $K_5$ and $K_{3,3}$.
Given a graph $G$, 
any subgraph of $G$ that is a subdivision of $K_5$ or $K_{3,3}$
is called a \emph{Kuratowski subgraph} of $G$.

Then one might ask if Kuratowski's result can be extended 
to higher genus surfaces in terms of Kuratowski subgraphs.
Glover has conjectured that if a finite graph $G$ is a minimal forbidden subgraph
for the nonorientable surface $\mathbb{N}_{\tilde{g}}$,
then $G$ can be written as the union of $\tilde{g}+1$ 
Kuratowski subgraphs such that 
the union of each pair of these fails to embed in the projective plane, 
the union of each triple of these fails to embed in the Klein bottle 
if $\tilde{g} \ge 2$, 
and the union of each triple of these fails to embed in the torus 
if $\tilde{g} \ge 3$.
It should be noted that $G$ is the union of 
$\tilde{g}+1$ Kuratowski subgraphs, i.e., every edge in $G$ is an edge 
in at least one of the Kuratowski subgraphs. 
The set of $\tilde{g}+1$ subgraphs described in the conjecture is called 
a \emph{Kuratowski covering} and the conjecture is called 
the \emph{Kuratowski covering conjecture}.  
In this paper, 
we prove the following restricted version of the above conjecture.

\begin{theorem}\label{th:main-theorem} The Kuratowski covering conjecture 
is true for every graph of order $<$ 10 
for $\mathbb{N}_{3}$ and $\mathbb{N}_{4}$. \end{theorem}

We prove this theorem by providing a Kuratowski covering for every minimal 
forbidden subgraph of order $<$ 10 for $\mathbb{N}_{3}$ and $\mathbb{N}_{4}$.

We use the complete lists of
minimal forbidden subgraphs of order $<$ 10 for some surfaces.
Archdeacon and Huneke \cite{AH} showed that
there are finitely many minimal forbidden subgraphs for nonorientable surfaces, 
and Robertson and Seymour \cite{RS} independently showed that 
there are finitely many minimal forbidden subgraphs for arbitrary surfaces.
Kuratowski \cite{K} showed that minimal forbidden subgraphs for the plane
are $K_5$ and $K_{3,3}$.
A list of minimal forbidden subgraphs for the projective plane
has been found by Glover, Huneke, and Wang \cite{GHW} and
Archdeacon \cite{A} proved that this list is complete.
However, for the higher genus surfaces, 
the complete lists of minimal forbidden subgraphs are not known.
The complete list of 8-vertex minimal forbidden subgraphs for the Klein bottle 
has been found by Huneke, McQuillan, and Richter \cite{HMR} and 
the complete list of 9-vertex minimal forbidden subgraphs for the Klein bottle
have been found by Cashy \cite{C} and Hur \cite{Hur1}.
The complete list of 8-vertex minimal forbidden subgraphs for the torus 
has been found by Duke and Haggard \cite{DH}, and 
the complete list of 9-vertex minimal forbidden subgraphs for the torus
has been found by Hlavacek \cite{H}.
We note that Wendy Myrvold in the Department of Computer Science, 
University of Victoria independently found about 235,000 minimal forbidden 
subgraphs for the torus using a computer \cite{M}, \cite{MN}.
In particular, the list of 9-vertex minimal forbidden subgraphs for the torus
found by Hlavacek coincides with the list of 9-vertex minimal forbidden subgraphs
for the torus found by Myrvold \cite{M}.

The remainder of this paper is organized as follows. 
In Section \ref{sec:n34}, we find Kuratowski coverings 
for all minimal forbidden subgraphs of order $<$ 10 
for $\mathbb{N}_{3}$ and $\mathbb{N}_{4}$, which prove
Theorem \ref{th:main-theorem}.

\medskip
{\it Remark 1.}
Theorem \ref{th:main-theorem} is part of the main result of author's thesis 
\cite{Hur1}, \cite{Hur3} : Every minimal forbidden subgraph of order $<$ 10 
satisfies the Kuratowski covering conjecture.

{\it Remark 2.}
A strengthened form of the Kuratowski covering conjecture analogous to the 
complete Kuratowski theorem for the plane says that 
a finite graph $G$ fails to embed in $\mathbb{N}_{\tilde{g}}$ 
if and only if there are $\tilde{g}+1$ Kuratowski subgraphs in $G$ 
satisfying the conditions of the Kuratowski covering conjecture.


\section{Kuratowski coverings for $\mathbb{N}_3$ and $\mathbb{N}_4$}
\label{sec:n34}

\subsection{Kuratowski coverings for $\mathbb{N}_3$}

\begin{figure}
\centering
\includegraphics[viewport=0 0 9in 14in, width=12cm]{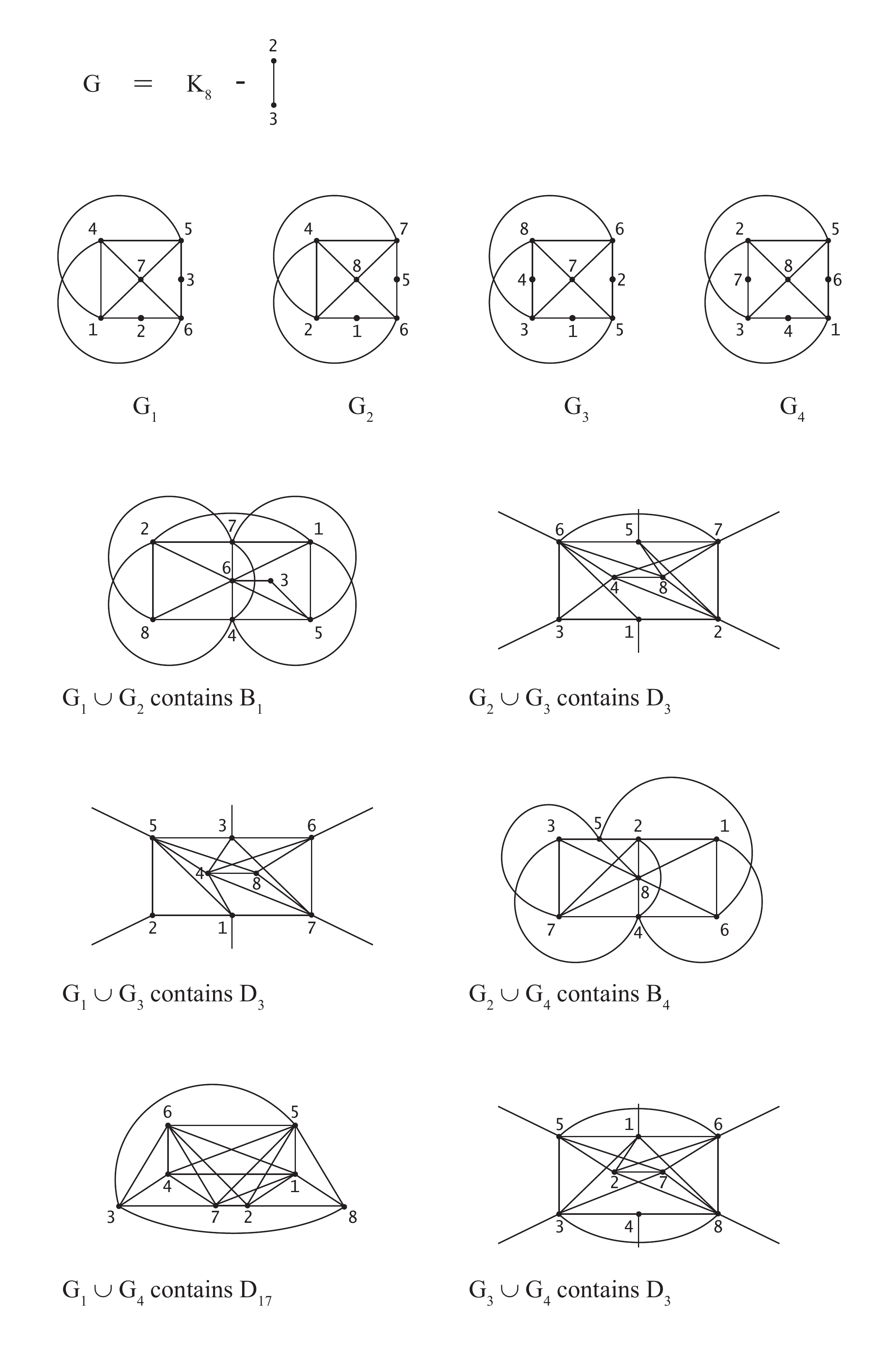}
\caption{\label{fig:kc-v8-n3-1-1} Kuratowski covering of $G = \tilde{I}^3_{8,1}$
for $\mathbb{N}_3$}
\end{figure}

\begin{figure}
\centering
\includegraphics[viewport=0 0 9in 13in, width=12cm]{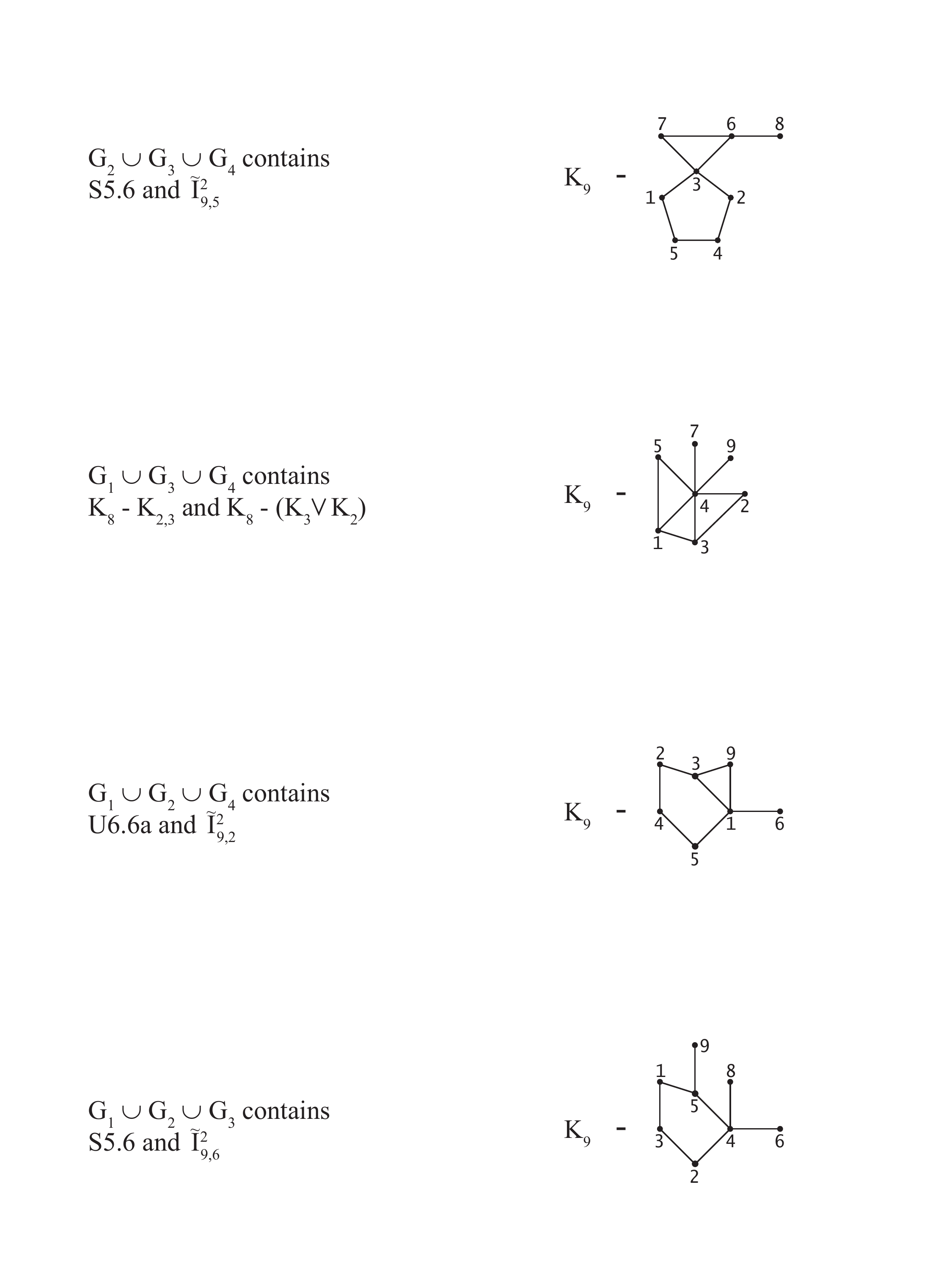}
\caption{\label{fig:kc-v9-n3-16-2} Kuratowski covering of $G = \tilde{I}^3_{9,16}$
for $\mathbb{N}_3$ (Continued)}
\end{figure}

The complete list of minimal forbidden subgraphs of order $<$ 10 
for $\mathbb{N}_{3}$ is given in \cite{Hur2}.
Since the genus of $\mathbb{N}_{3}$ is three, for each minimal forbidden 
subgraph $G$ of order $<$ 10 for $\mathbb{N}_{3}$, 
we find four Kuratowski subgraphs $ G_1 , G_2 , G_3 , G_4 $ 
as a Kuratowski covering such that 
the union of every pair of these contains a subdivision of 
a minimal forbidden subgraph for the projective plane and
the union of every triple of these contains a subdivision of 
a minimal forbidden subgraph for the torus 
and a subdivision of a minimal forbidden subgraph for the Klein bottle.
We show these Kuratowski coverings in Figure \ref{fig:kc-v8-n3-1-1}, $\ldots$ ,
\ref{fig:kc-v9-n3-16-2}.
The names of minimal forbidden subgraphs for the projective are from \cite{GHW}.
The 8-vertex minimal forbidden subgraphs for the torus are
$K_8-K_3$, $K_8 - (K_{1,2} \cup 2K_2)$, and $K_8 - K_{2,3}$ \cite{DH},
and the 8-vertex minimal forbidden subgraphs for the Klein bottle are
$K_8 - 4K_2$, $K_8 - (K_3 \vee K_2)$, $K_8 - 2K_3$,
$K_8 - 2K_{1,3}$, and $K_8 - (K_{1,4} \cup K_3)$,
which are from \cite{HMR}.
The names of 9-vertex minimal forbidden subgraphs for the torus are 
from \cite{H} and there are 63 9-vertex minimal forbidden subgraphs 
$\tilde{I}^2_{9,1}, \ldots , \tilde{I}^2_{9,63}$ for the Klein bottle 
\cite{C}, \cite{Hur1}.

\subsection{Kuratowski coverings for $\mathbb{N}_4$}

\begin{figure}[h]
\centering
\includegraphics[viewport=0 0 8.5in 14in, width=11cm]{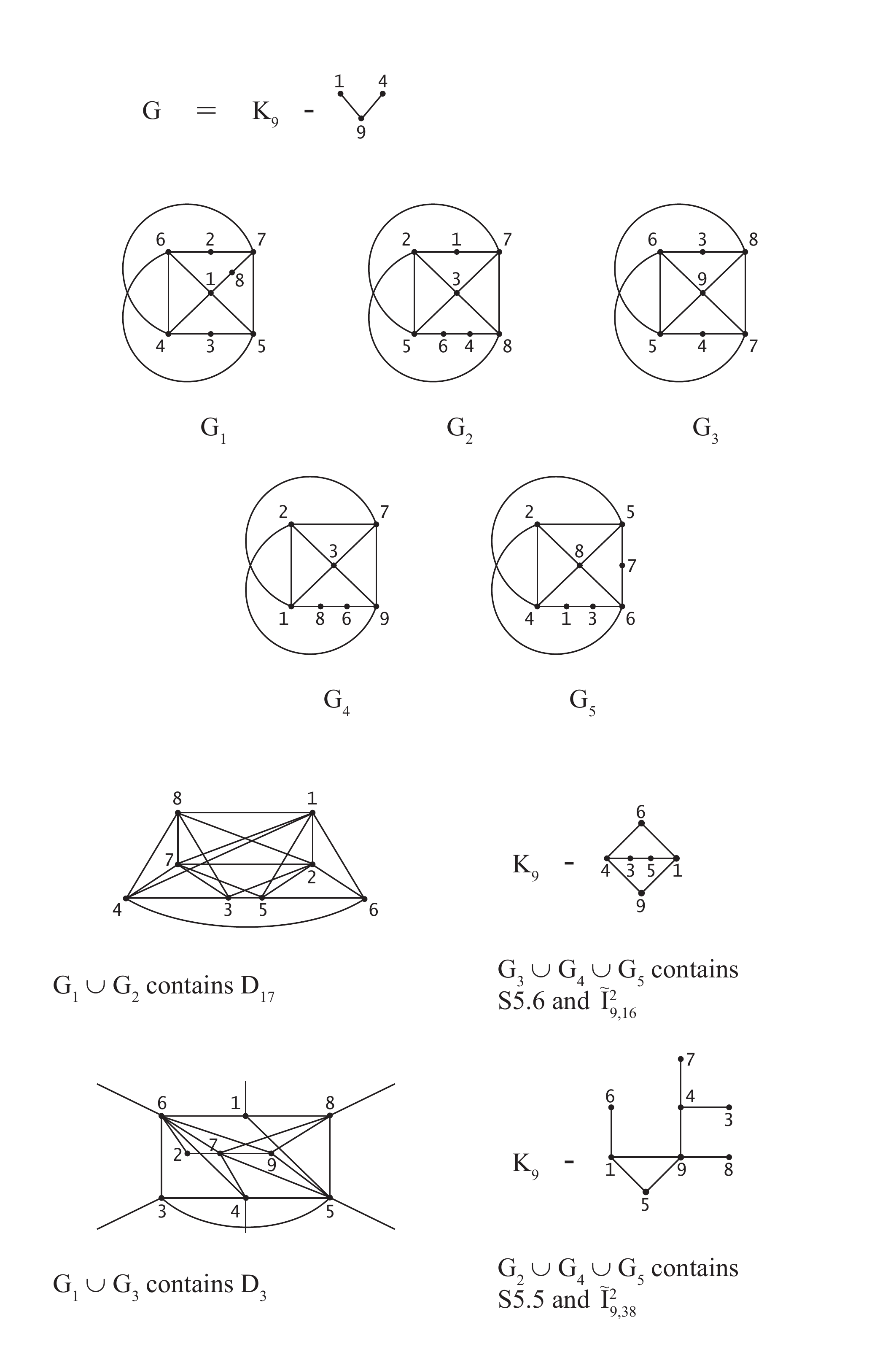}
\caption{\label{fig:kc-v9-n4-1-1} Kuratowski covering of $G = \tilde{I}^4_{9,1}$
for $\mathbb{N}_4$}
\end{figure}

\begin{figure}
\centering
\includegraphics[viewport=0 0 8.5in 14in, width=11cm]{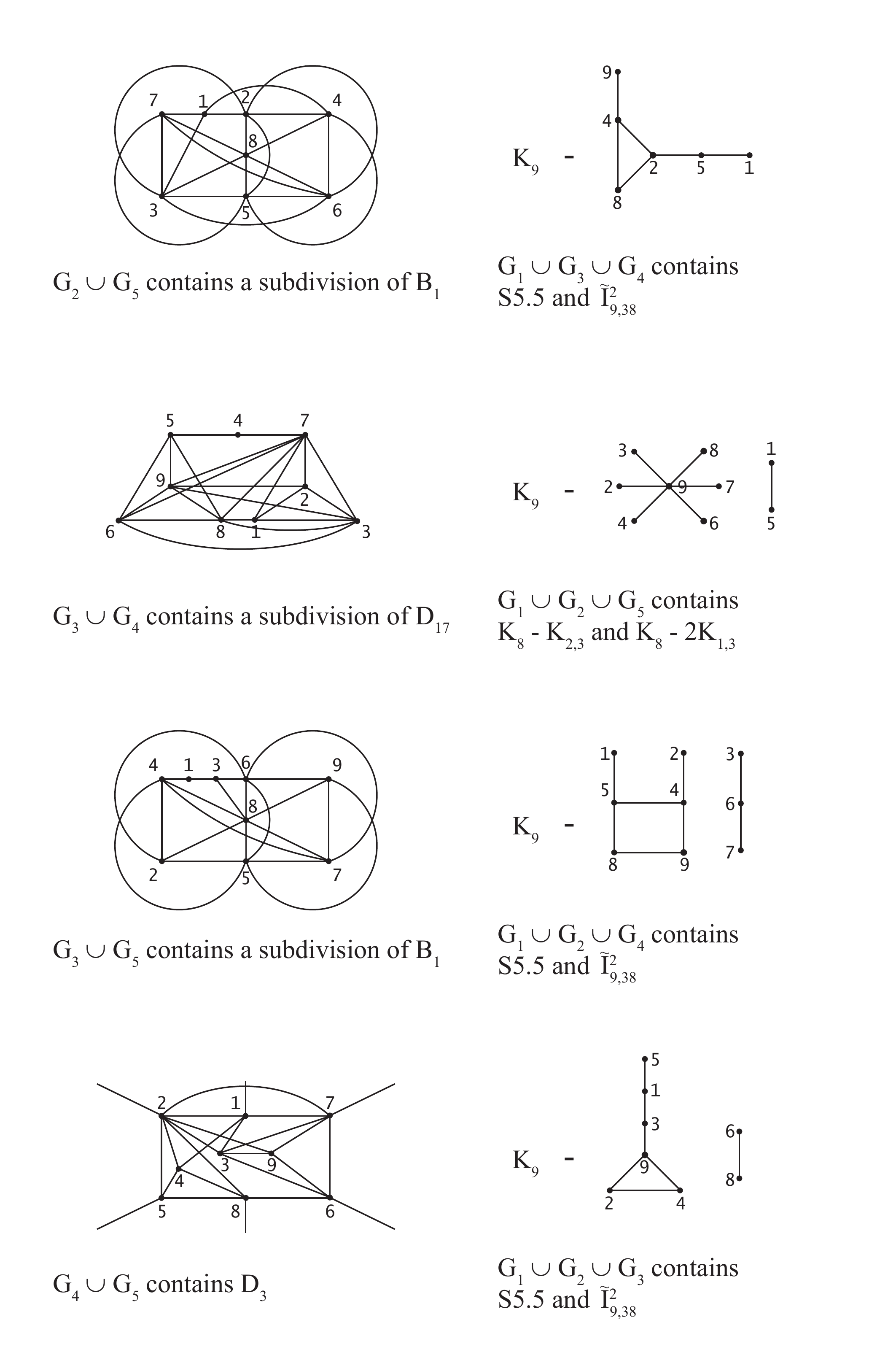}
\caption{\label{fig:kc-v9-n4-2-3} Kuratowski covering of $G = \tilde{I}^4_{9,2}$
for $\mathbb{N}_4$ (Continued)}
\end{figure}

There are two minimal forbidden subgraphs $K_9 - K_{1,2}$ and $K_9 - 2K_2$
of order $<$ 10 for $\mathbb{N}_4$ \cite{Hur2} and we need five 
Kuratowski subgraphs $ G_1 , G_2 , G_3 , G_4 , G_5  $ 
as a Kuratowski covering which satisfies 
the conditions of the Kuratowski covering conjecture.
We show Kuratowski coverings of $K_9 - K_{1,2}$ and $K_9 - 2K_2$ 
in Figure \ref{fig:kc-v9-n4-1-1}, $\ldots$ , \ref{fig:kc-v9-n4-2-3}.

\begin{figure}
\centering
\includegraphics[viewport=0 0 9in 13in, width=12cm]{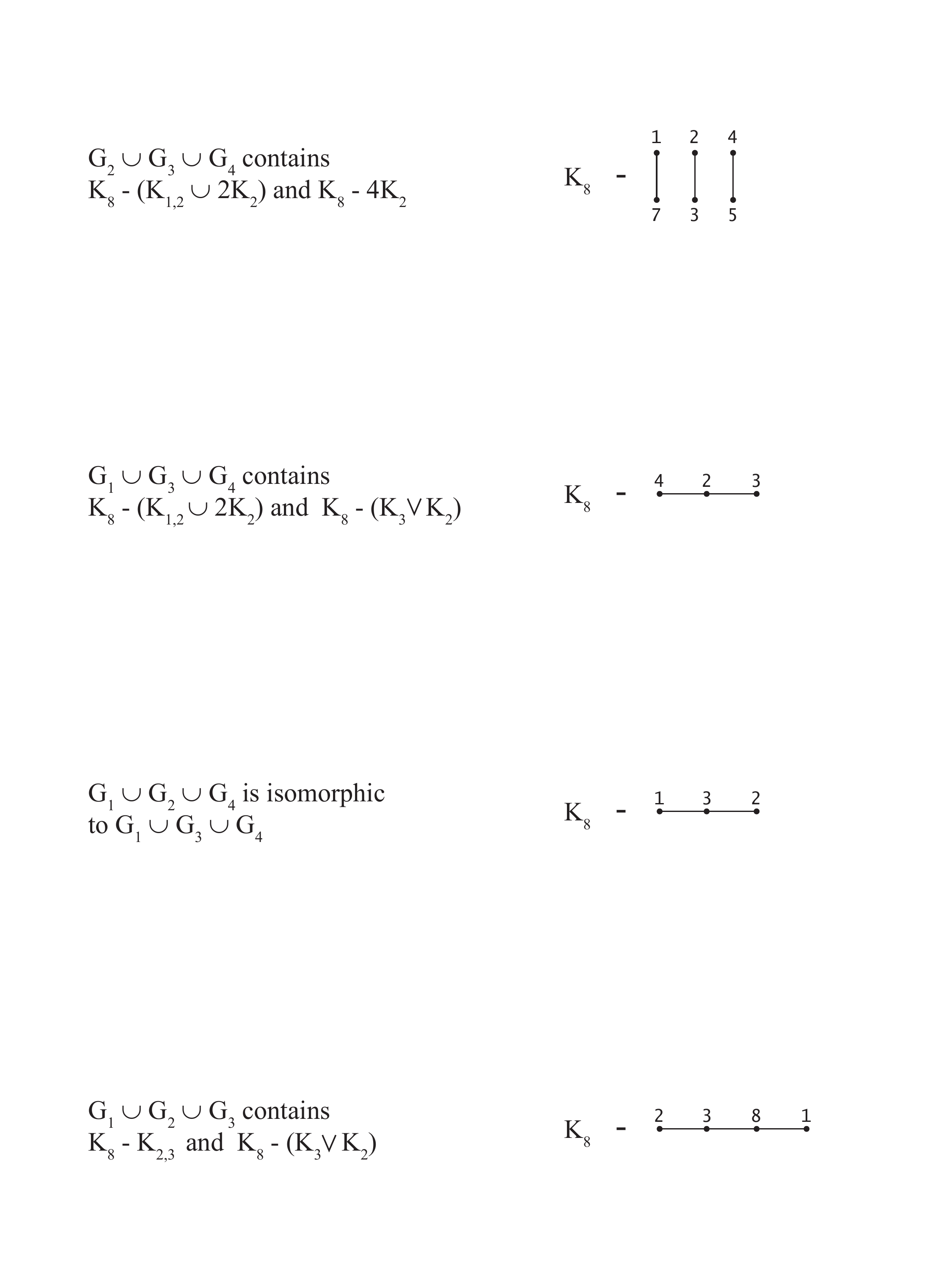}
\caption{\label{fig:kc-v8-n3-1-2} Kuratowski covering of $G = \tilde{I}^3_{8,1}$
for $\mathbb{N}_3$ (Continued)}
\end{figure}

\begin{figure}
\centering
\includegraphics[viewport=0 0 9in 14in, width=12cm]{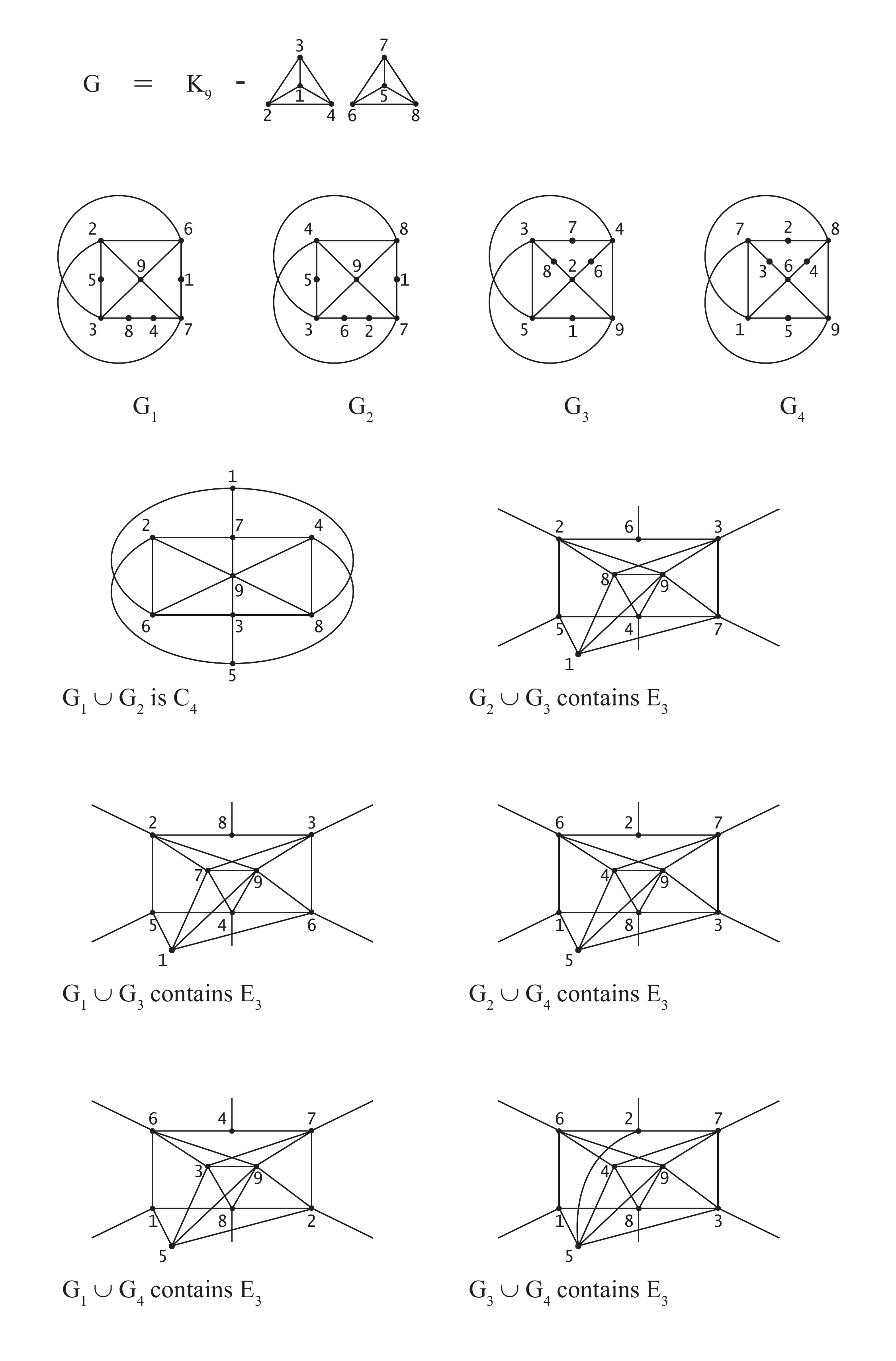}
\caption{\label{fig:kc-v9-n3-1-1} Kuratowski covering of $G = \tilde{I}^3_{9,1}$
for $\mathbb{N}_3$}
\end{figure}

\begin{figure}
\centering
\includegraphics[viewport=0 0 9in 13in, width=12cm]{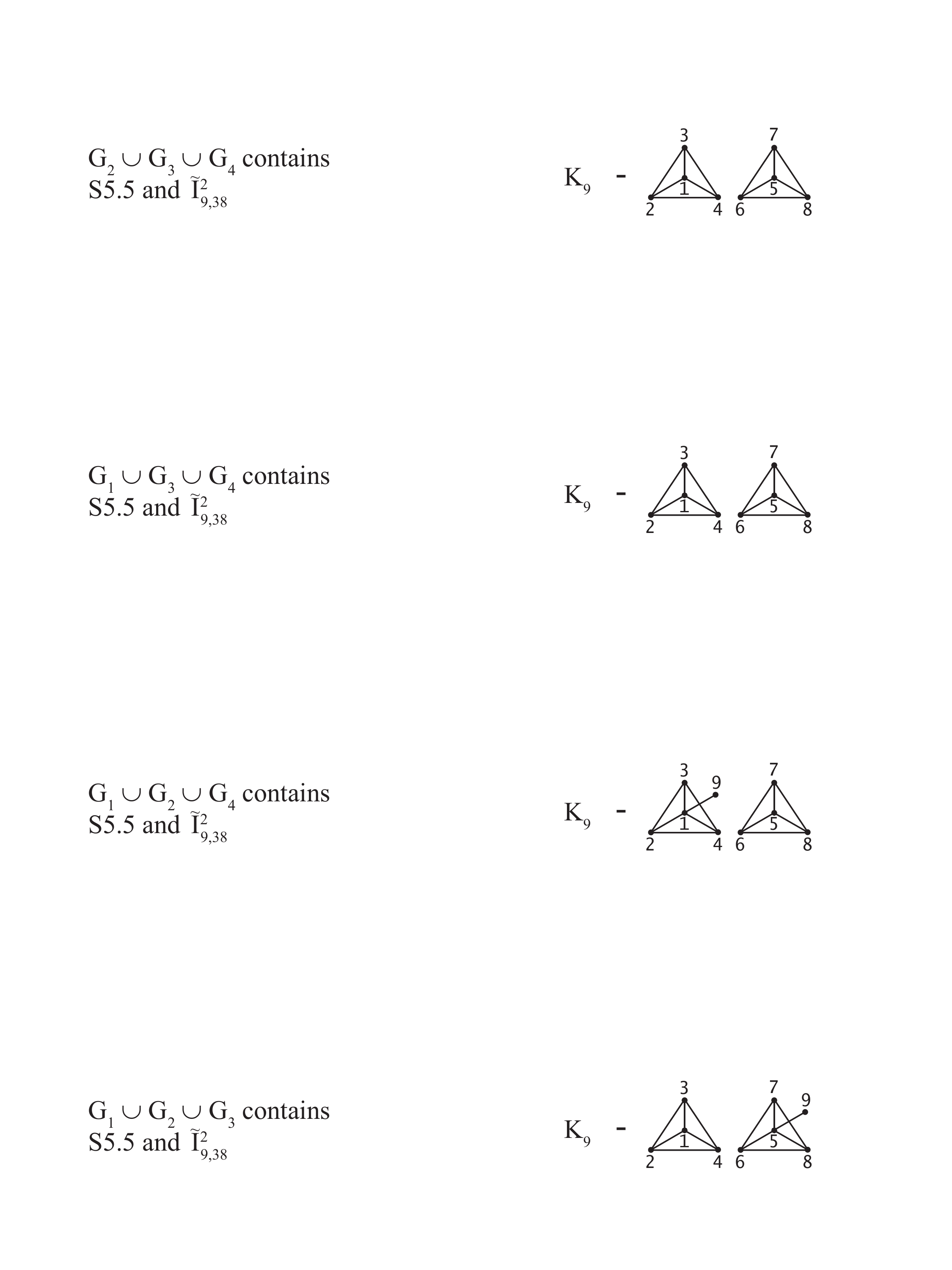}
\caption{\label{fig:kc-v9-n3-1-2} Kuratowski covering of $G = \tilde{I}^3_{9,1}$
for $\mathbb{N}_3$ (Continued)}
\end{figure}

\begin{figure}
\centering
\includegraphics[viewport=0 0 9in 14in, width=12cm]{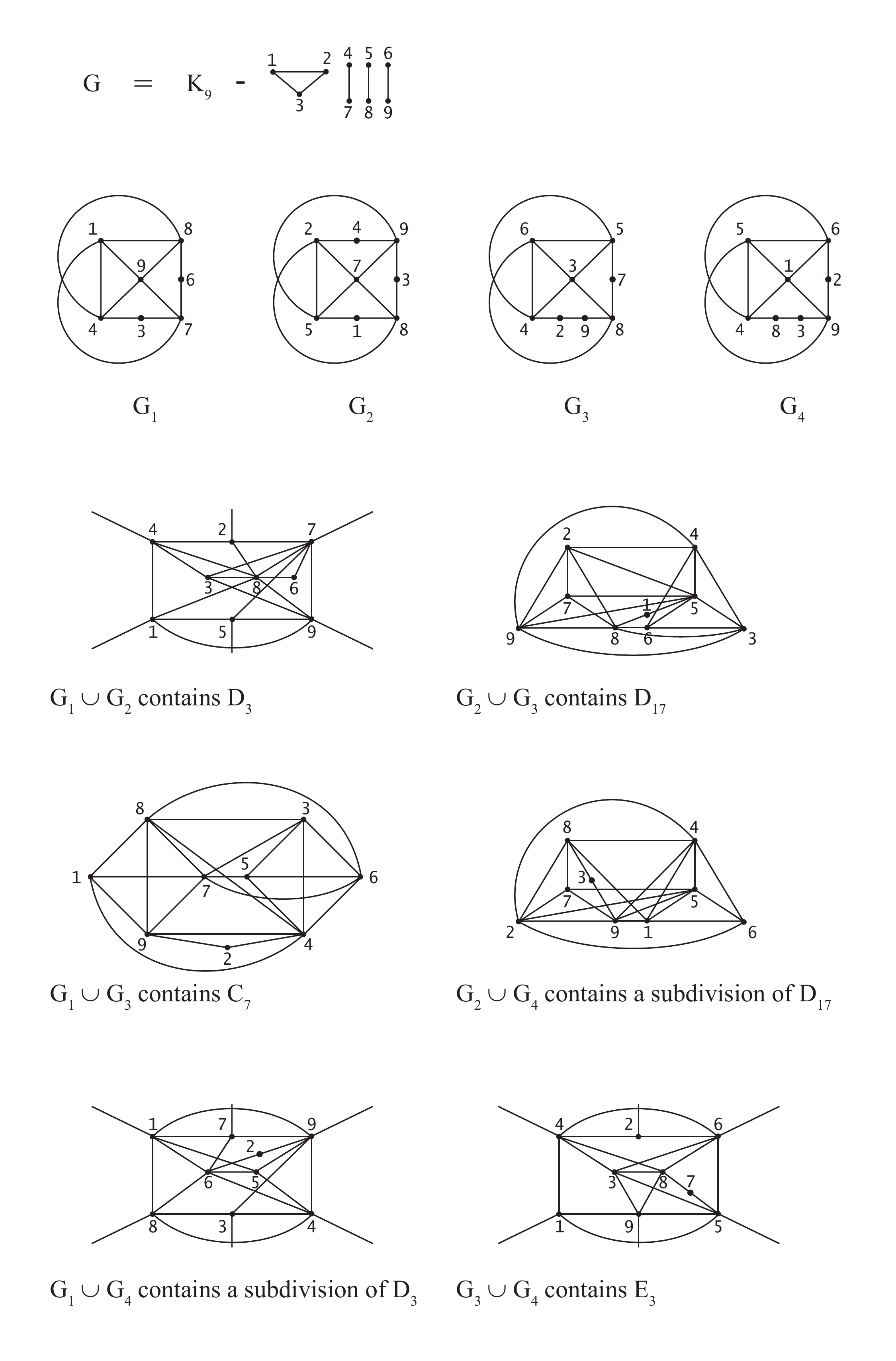}
\caption{\label{fig:kc-v9-n3-2-1} Kuratowski covering of $G = \tilde{I}^3_{9,2}$
for $\mathbb{N}_3$}
\end{figure}

\begin{figure}
\centering
\includegraphics[viewport=0 0 9in 13in, width=12cm]{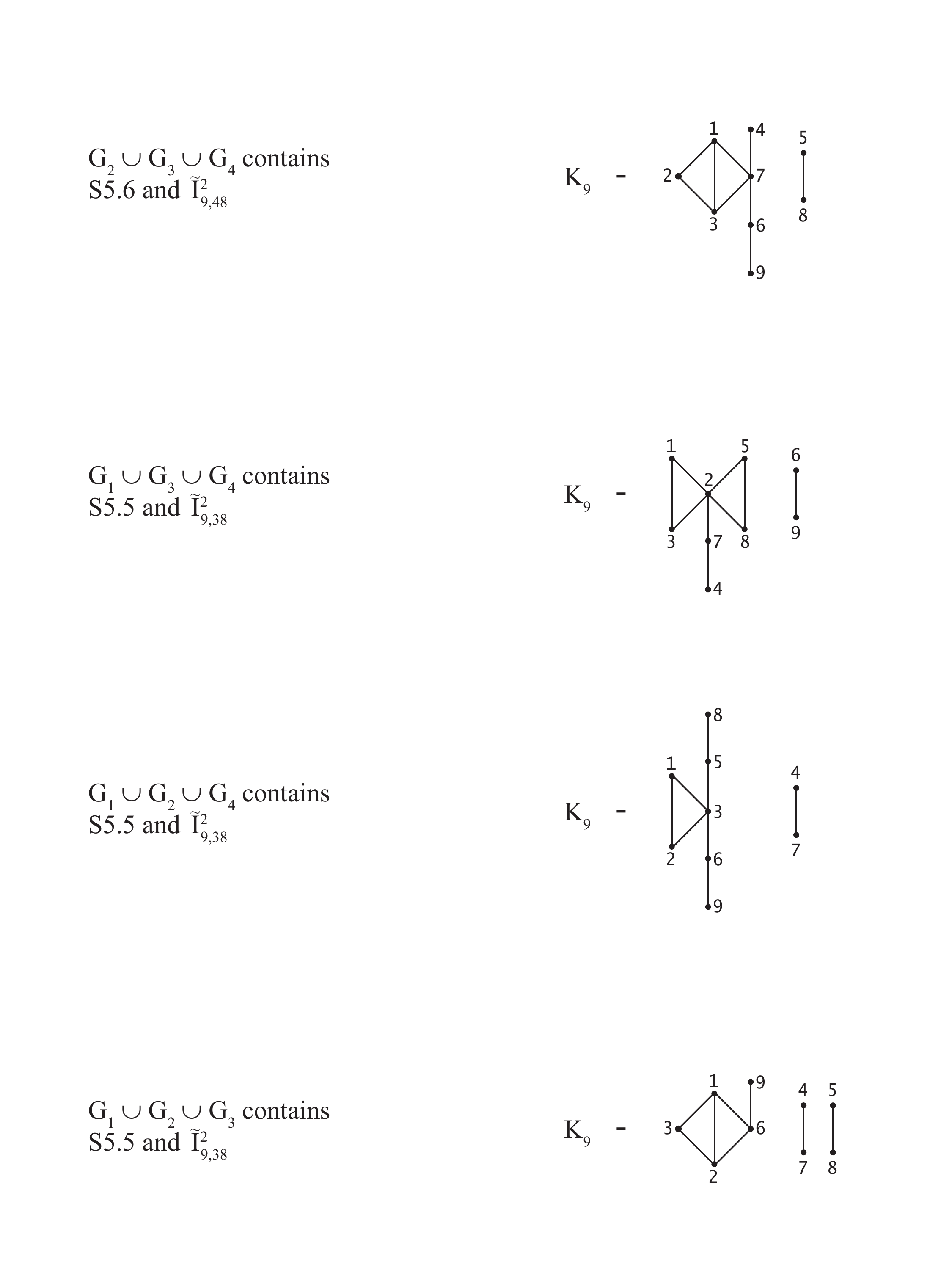}
\caption{\label{fig:kc-v9-n3-2-2} Kuratowski covering of $G = \tilde{I}^3_{9,2}$
for $\mathbb{N}_3$ (Continued)}
\end{figure}

\begin{figure}
\centering
\includegraphics[viewport=0 0 9in 14in, width=12cm]{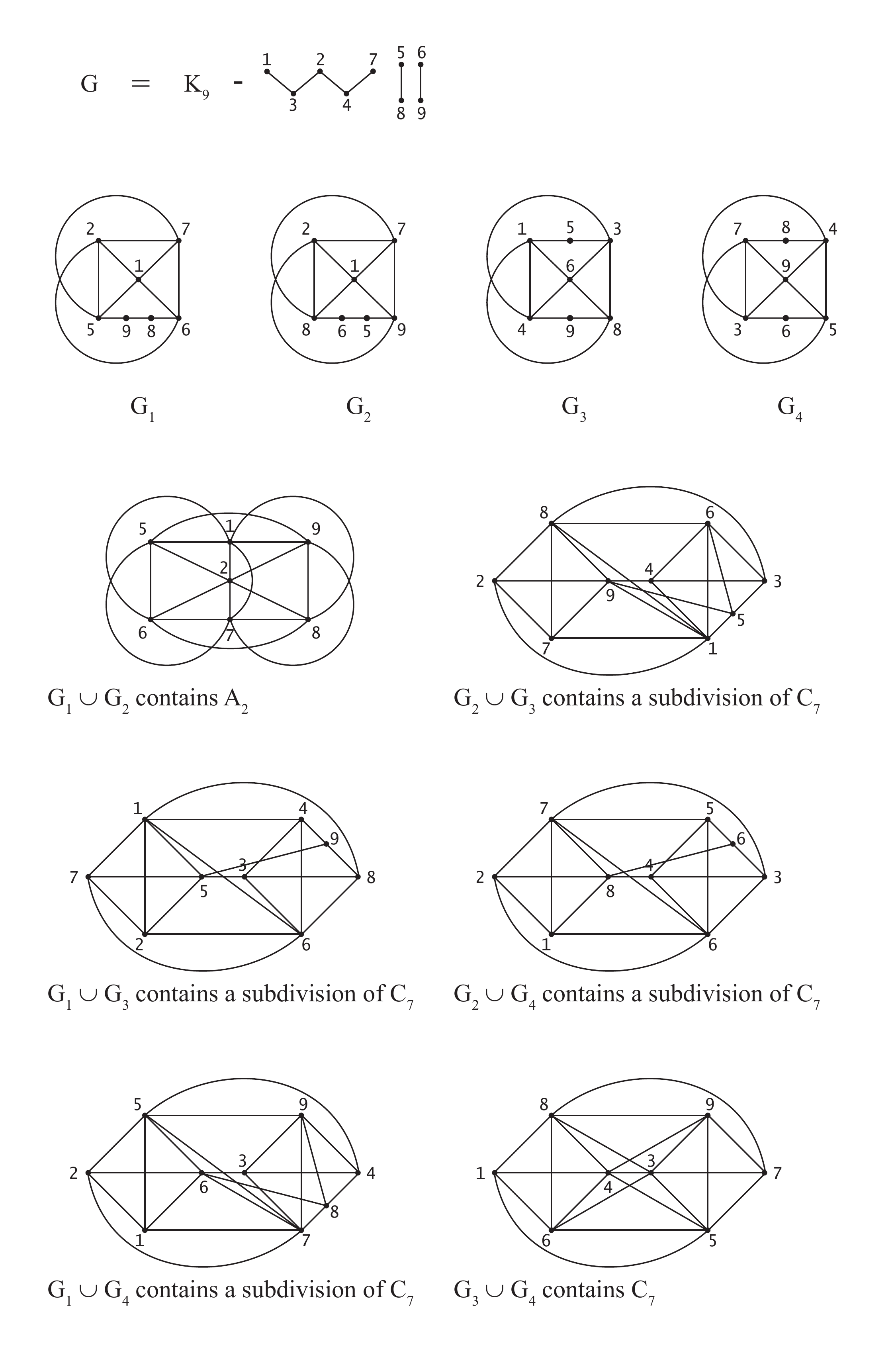}
\caption{\label{fig:kc-v9-n3-3-1} Kuratowski covering of $G = \tilde{I}^3_{9,3}$
for $\mathbb{N}_3$}
\end{figure}

\begin{figure}
\centering
\includegraphics[viewport=0 0 9in 13in, width=12cm]{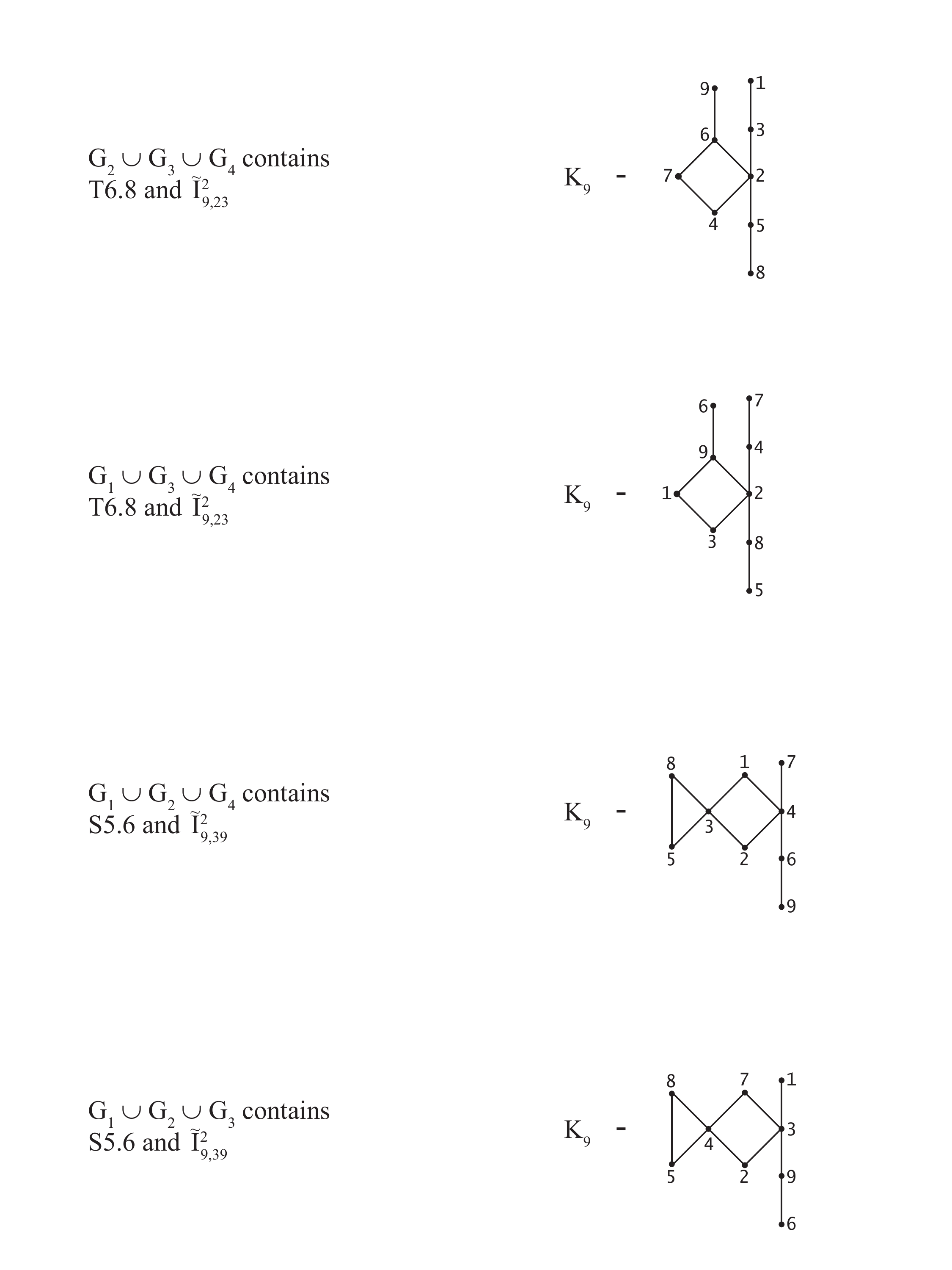}
\caption{\label{fig:kc-v9-n3-3-2} Kuratowski covering of $G = \tilde{I}^3_{9,3}$
for $\mathbb{N}_3$ (Continued)}
\end{figure}

\begin{figure}
\centering
\includegraphics[viewport=0 0 9in 14in, width=12cm]{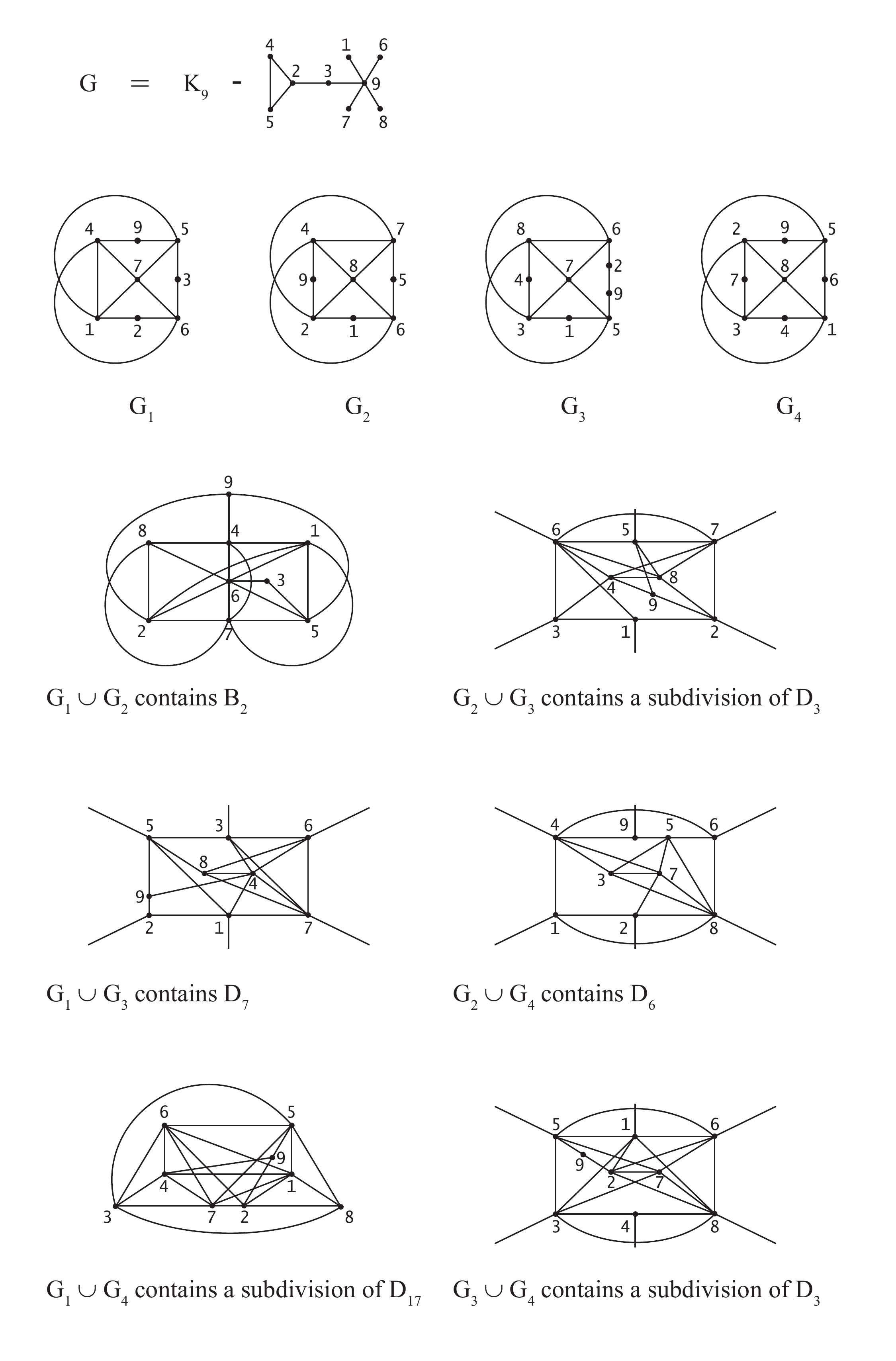}
\caption{\label{fig:kc-v9-n3-4-1} Kuratowski covering of $G = \tilde{I}^3_{9,4}$
for $\mathbb{N}_3$}
\end{figure}

\begin{figure}
\centering
\includegraphics[viewport=0 0 9in 13in, width=12cm]{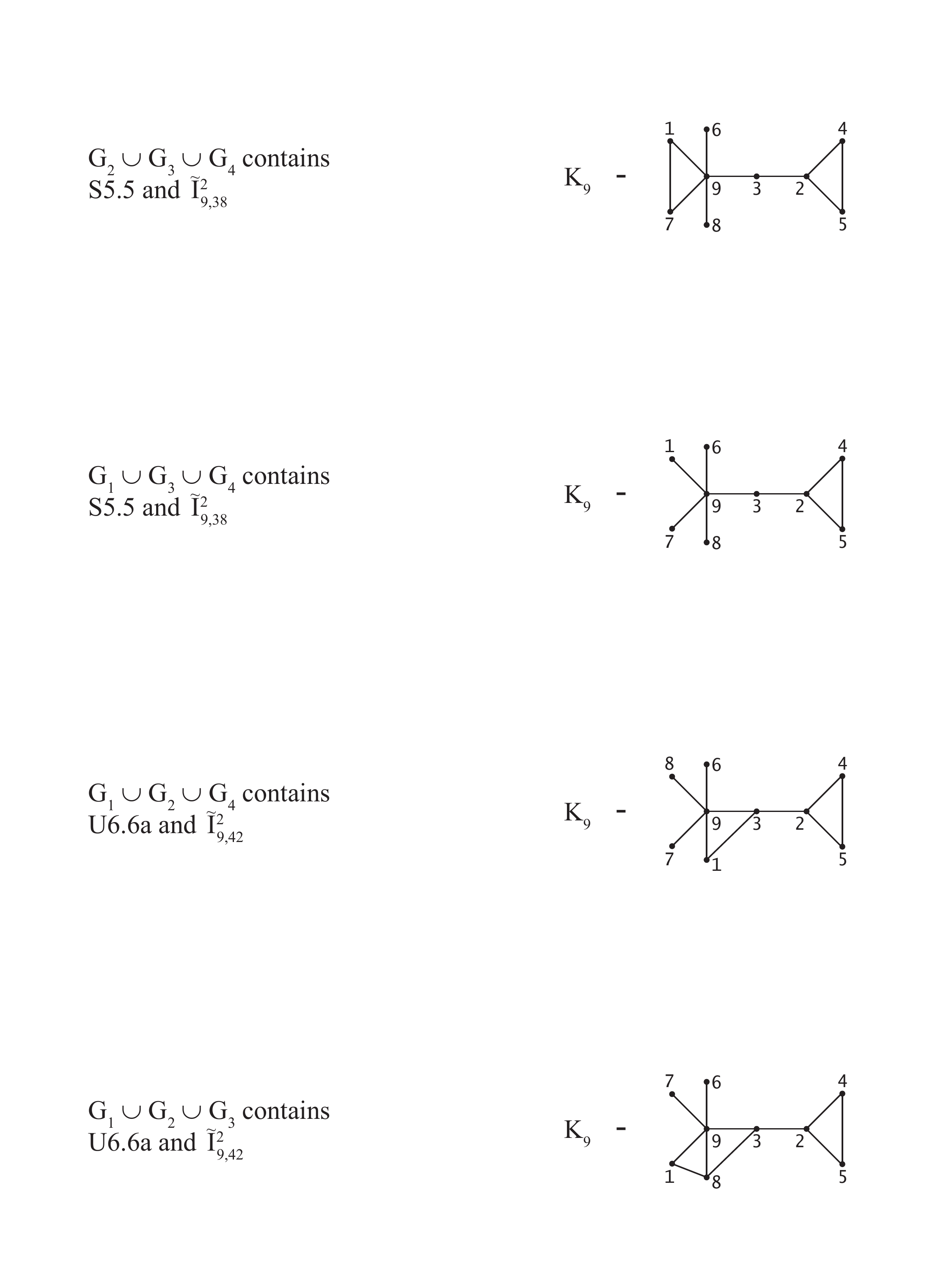}
\caption{\label{fig:kc-v9-n3-4-2} Kuratowski covering of $G = \tilde{I}^3_{9,4}$
for $\mathbb{N}_3$ (Continued)}
\end{figure}

\begin{figure}
\centering
\includegraphics[viewport=0 0 9in 14in, width=12cm]{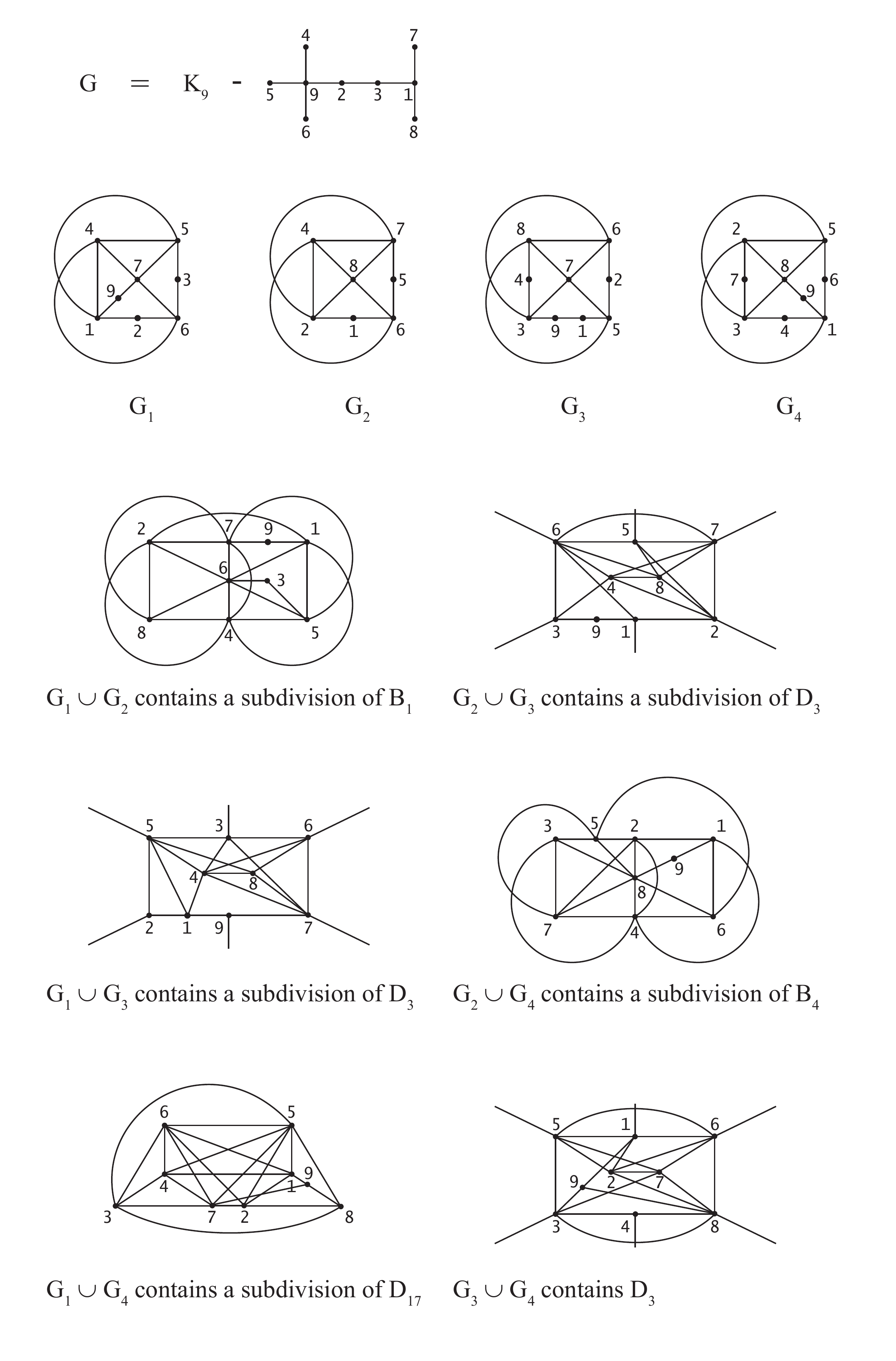}
\caption{\label{fig:kc-v9-n3-5-1} Kuratowski covering of $G = \tilde{I}^3_{9,5}$
for $\mathbb{N}_3$}
\end{figure}

\begin{figure}
\centering
\includegraphics[viewport=0 0 9in 13in, width=12cm]{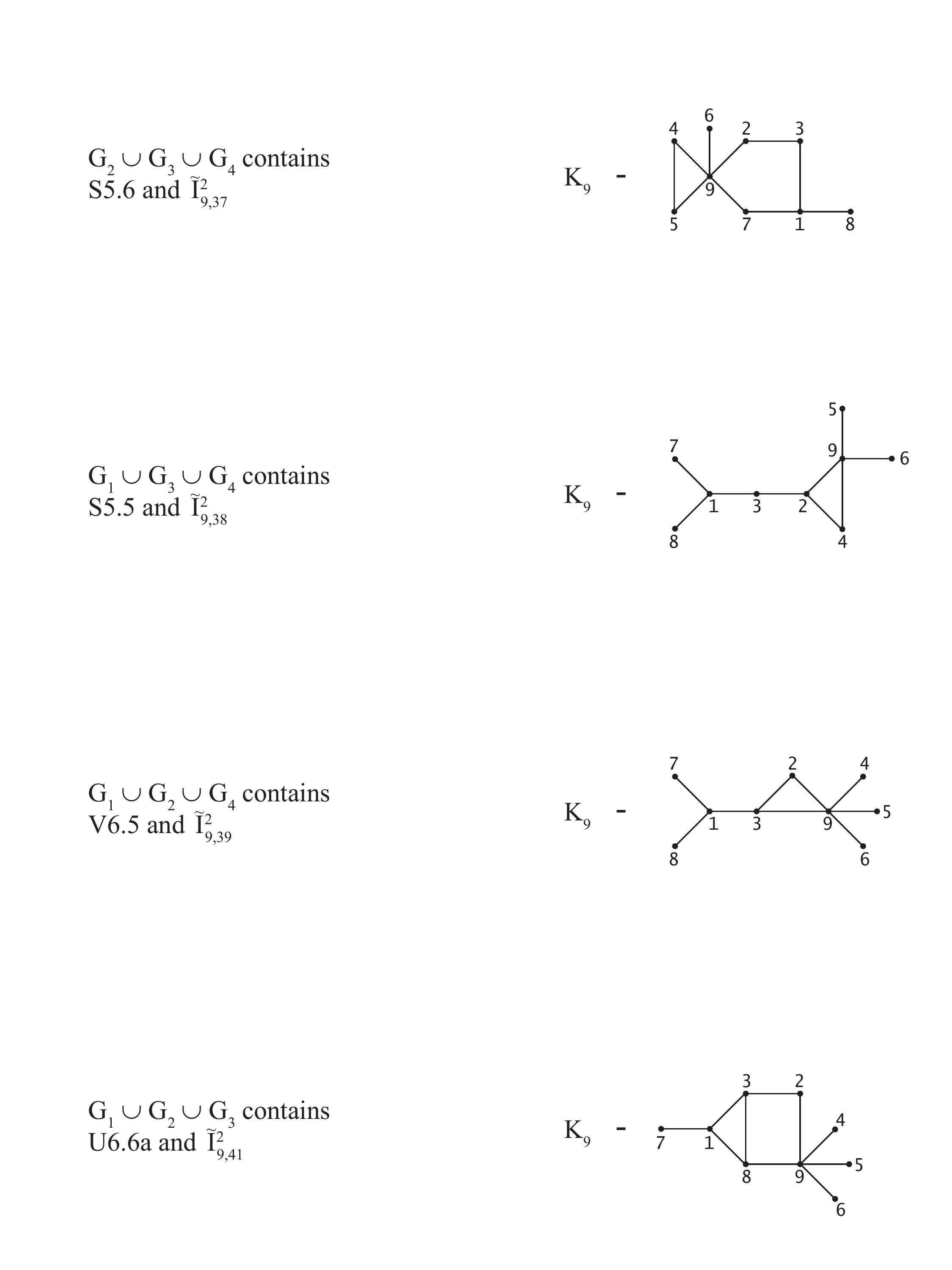}
\caption{\label{fig:kc-v9-n3-5-2} Kuratowski covering of $G = \tilde{I}^3_{9,5}$
for $\mathbb{N}_3$ (Continued)}
\end{figure}

\begin{figure}
\centering
\includegraphics[viewport=0 0 9in 14in, width=12cm]{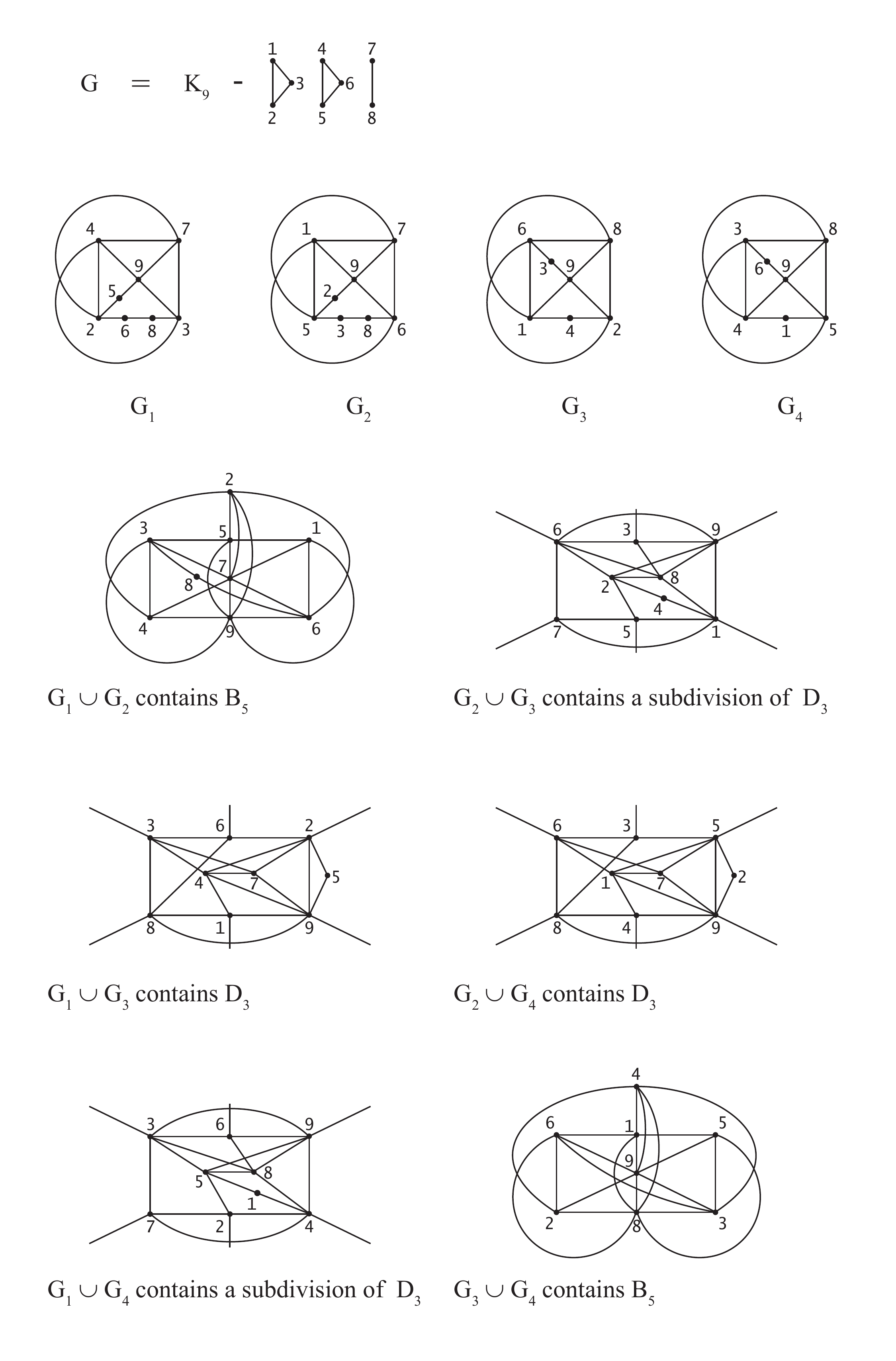}
\caption{\label{fig:kc-v9-n3-6-1} Kuratowski covering of $G = \tilde{I}^3_{9,6}$
for $\mathbb{N}_3$}
\end{figure}

\begin{figure}
\centering
\includegraphics[viewport=0 0 9in 13in, width=12cm]{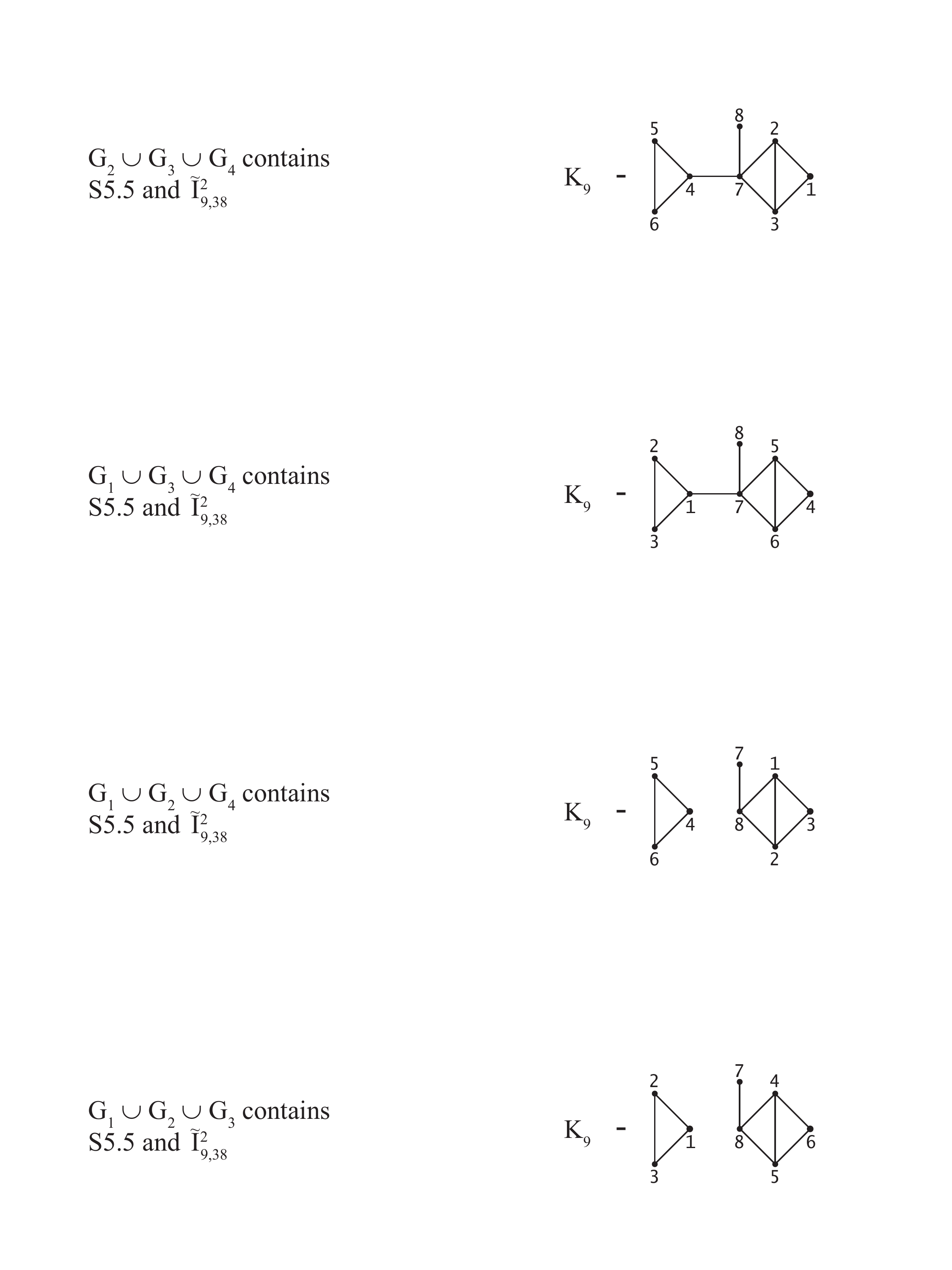}
\caption{\label{fig:kc-v9-n3-6-2} Kuratowski covering of $G = \tilde{I}^3_{9,6}$
for $\mathbb{N}_3$ (Continued)}
\end{figure}

\begin{figure}
\centering
\includegraphics[viewport=0 0 9in 14in, width=12cm]{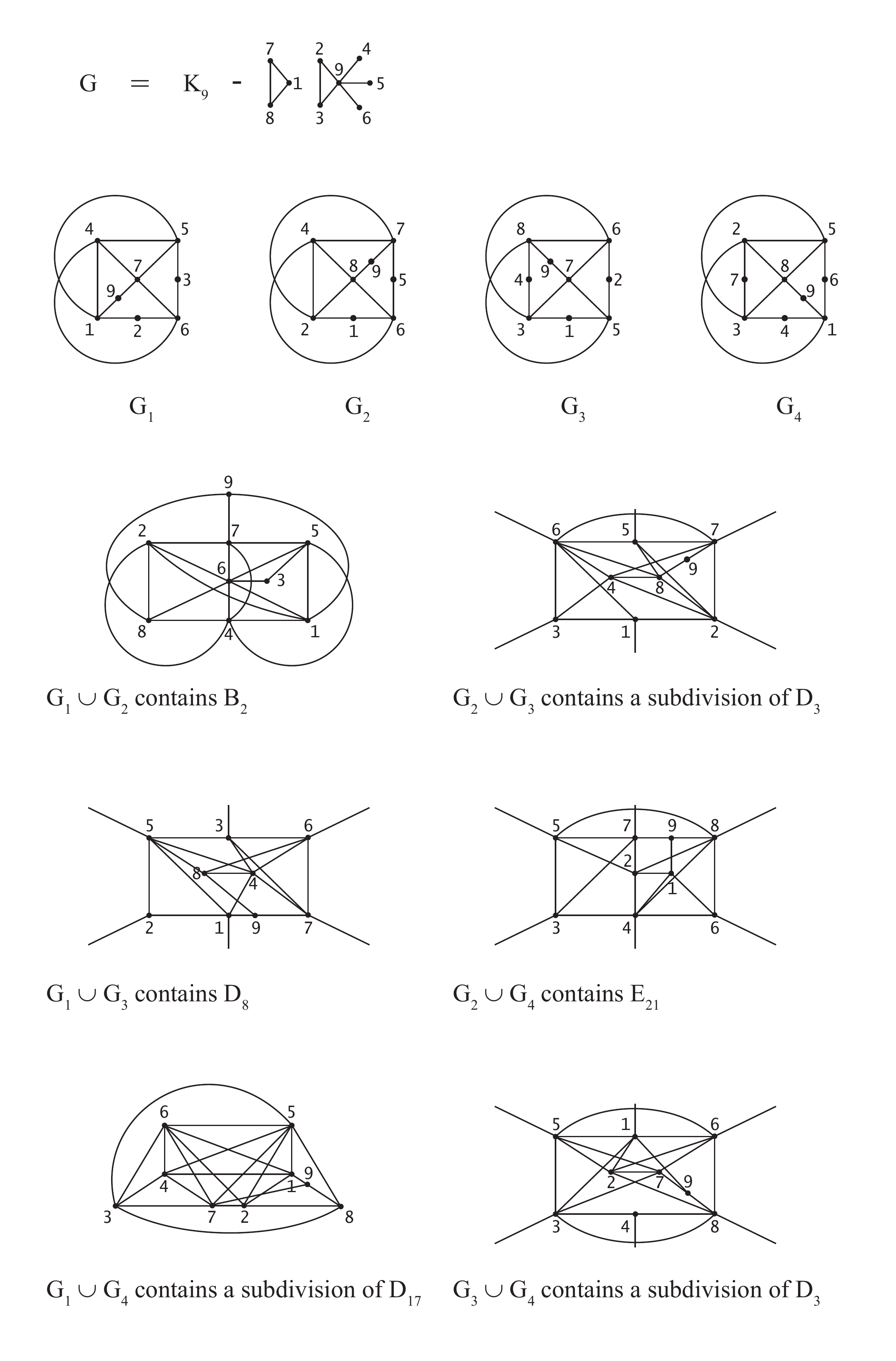}
\caption{\label{fig:kc-v9-n3-7-1} Kuratowski covering of $G = \tilde{I}^3_{9,7}$
for $\mathbb{N}_3$}
\end{figure}

\begin{figure}
\centering
\includegraphics[viewport=0 0 9in 13in, width=12cm]{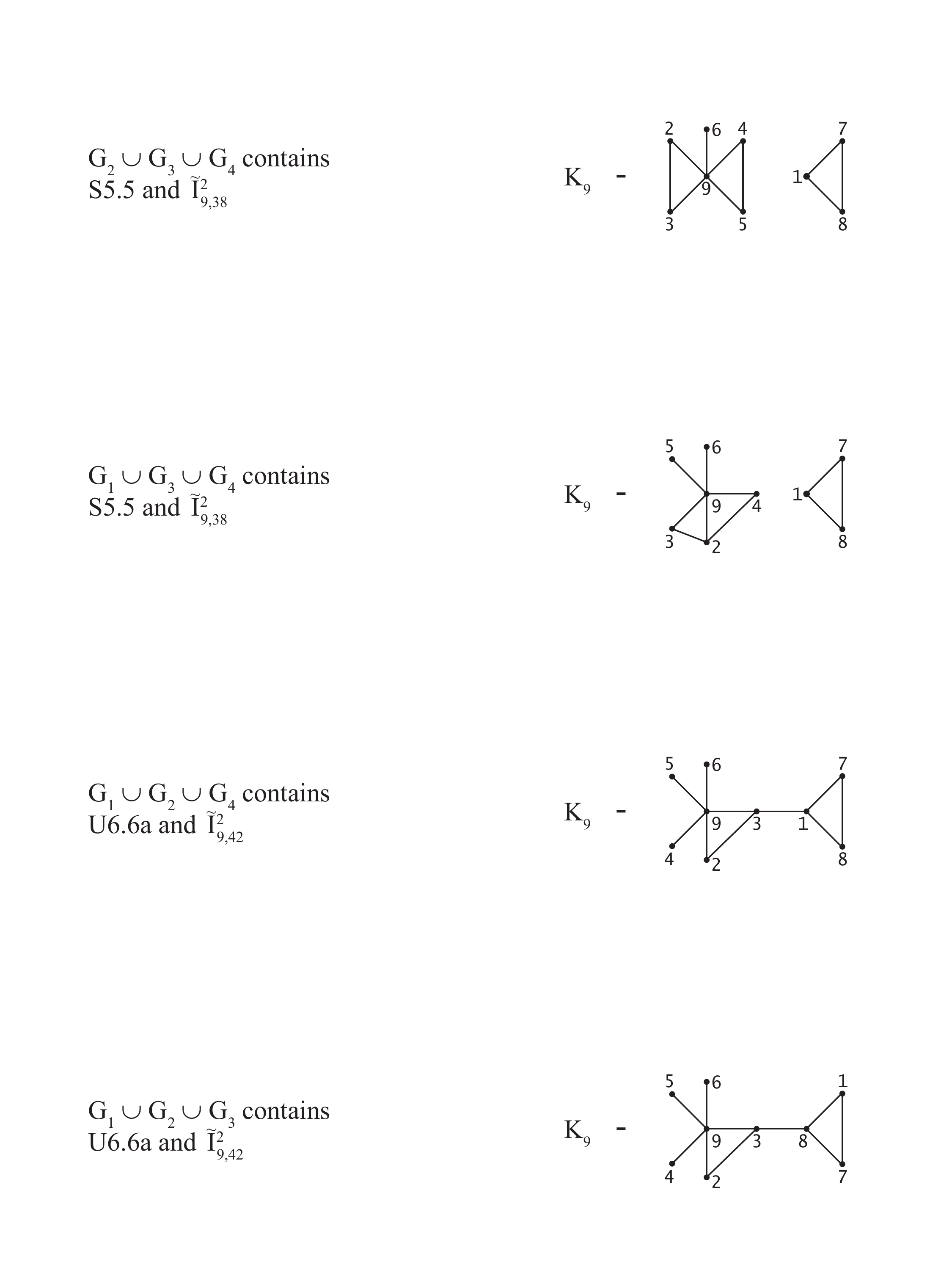}
\caption{\label{fig:kc-v9-n3-7-2} Kuratowski covering of $G = \tilde{I}^3_{9,7}$
for $\mathbb{N}_3$ (Continued)}
\end{figure}

\begin{figure}
\centering
\includegraphics[viewport=0 0 9in 14in, width=12cm]{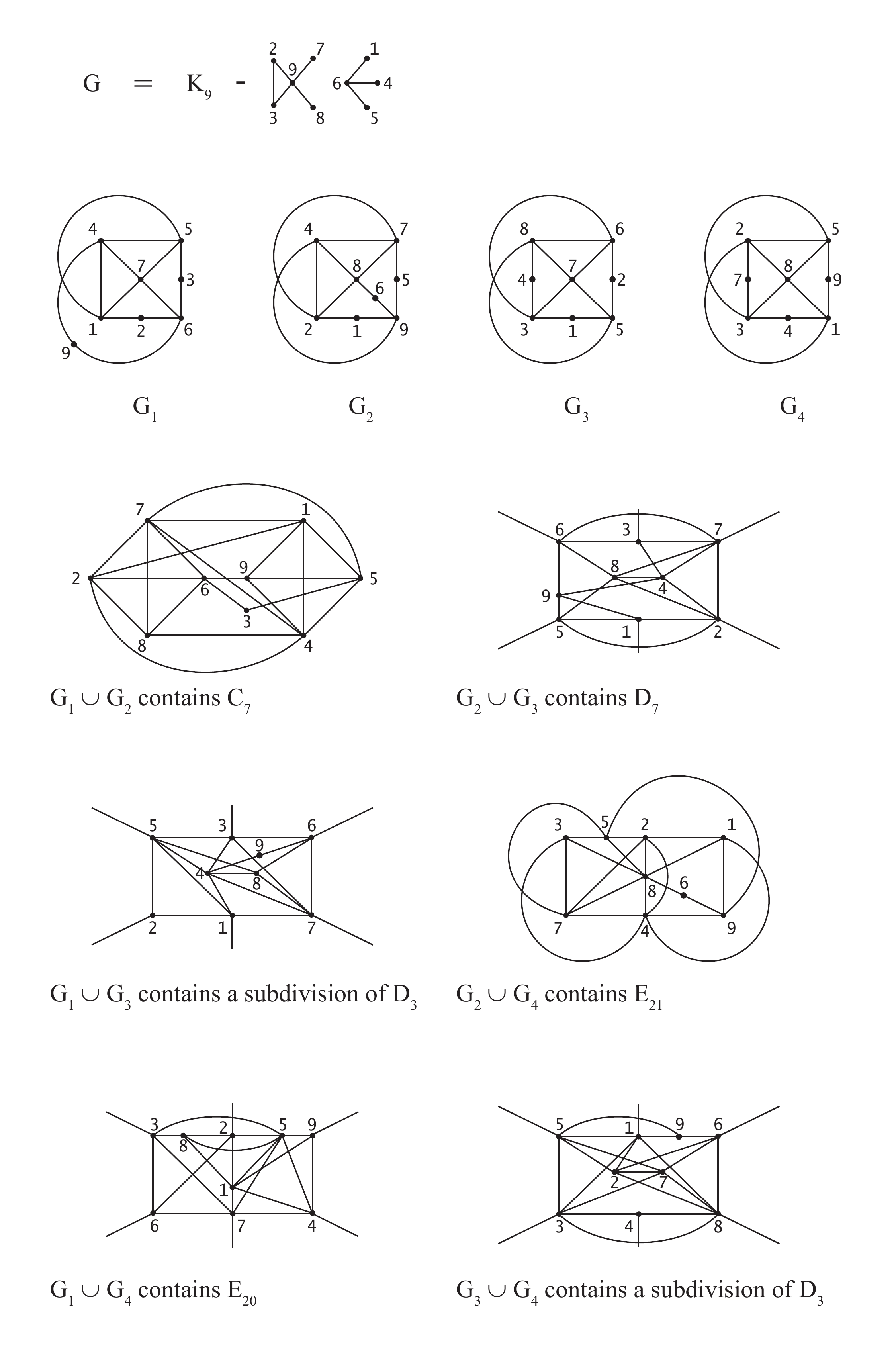}
\caption{\label{fig:kc-v9-n3-8-1} Kuratowski covering of $G = \tilde{I}^3_{9,8}$
for $\mathbb{N}_3$}
\end{figure}

\begin{figure}
\centering
\includegraphics[viewport=0 0 9in 13in, width=12cm]{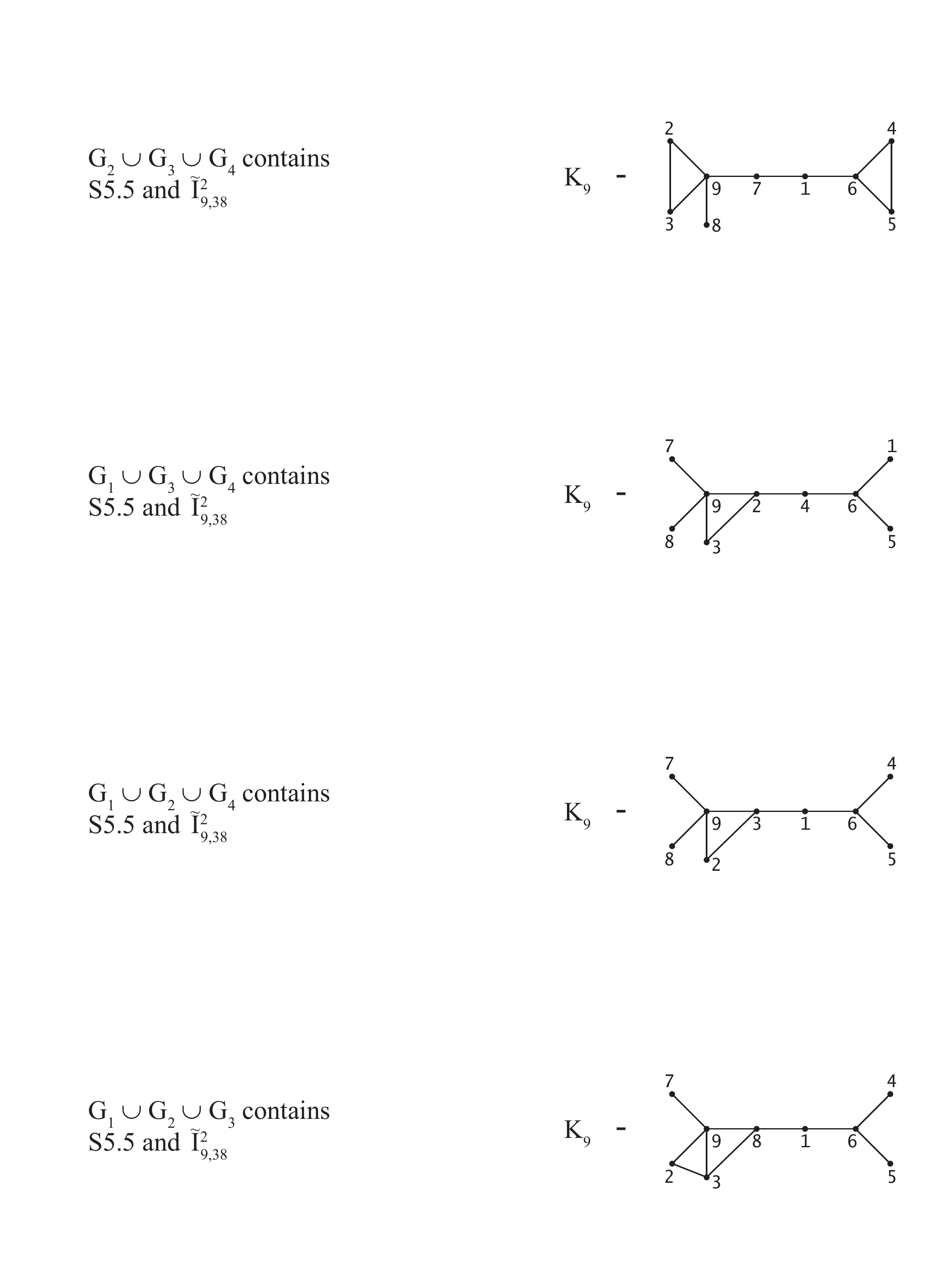}
\caption{\label{fig:kc-v9-n3-8-2} Kuratowski covering of $G = \tilde{I}^3_{9,8}$
for $\mathbb{N}_3$ (Continued)}
\end{figure}

\clearpage

\begin{figure}
\centering
\includegraphics[viewport=0 0 9in 14in, width=12cm]{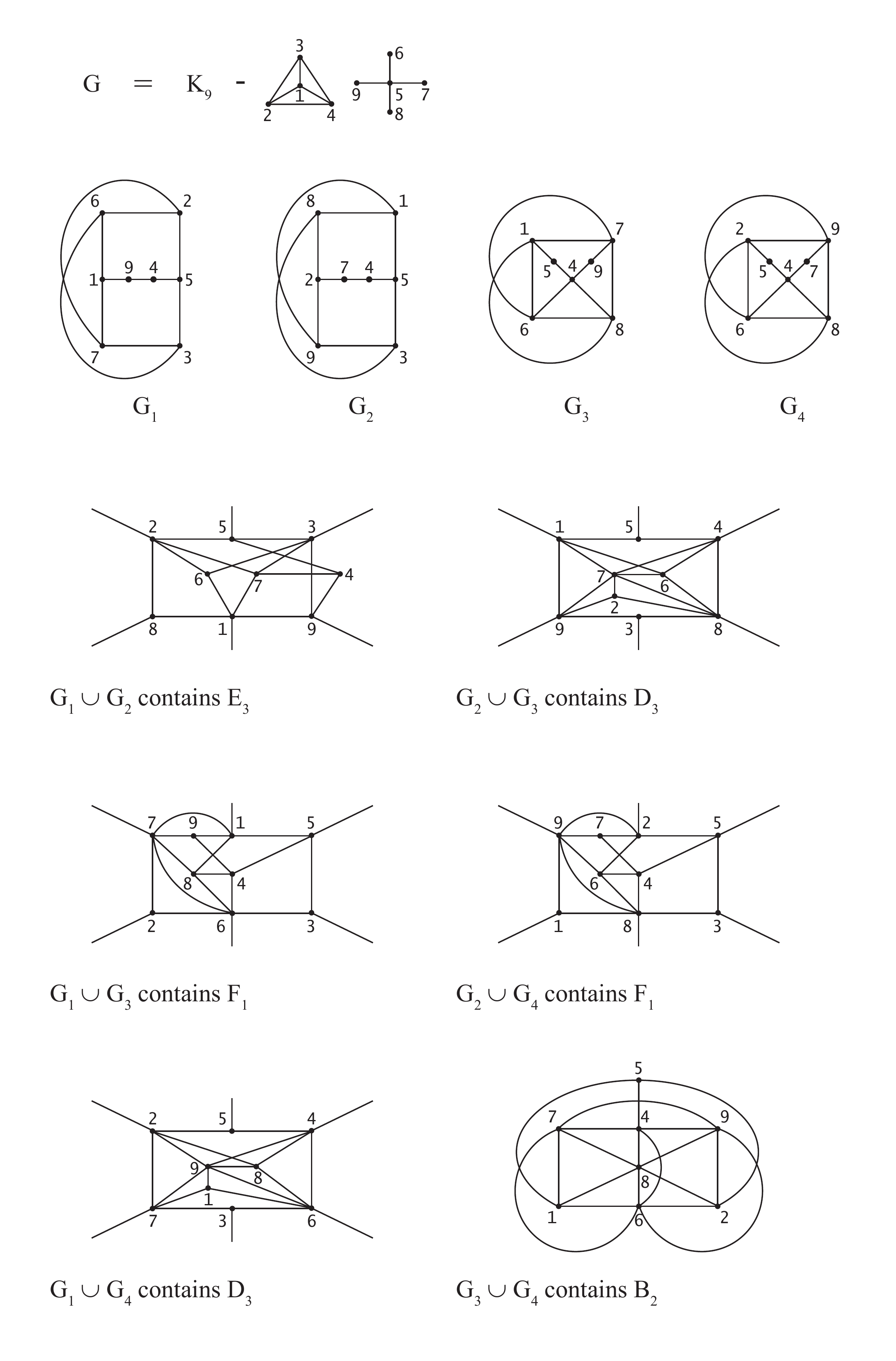}
\caption{\label{fig:kc-v9-n3-9-1} Kuratowski covering of $G = \tilde{I}^3_{9,9}$
for $\mathbb{N}_3$}
\end{figure}

\begin{figure}
\centering
\includegraphics[viewport=0 0 9in 13in, width=12cm]{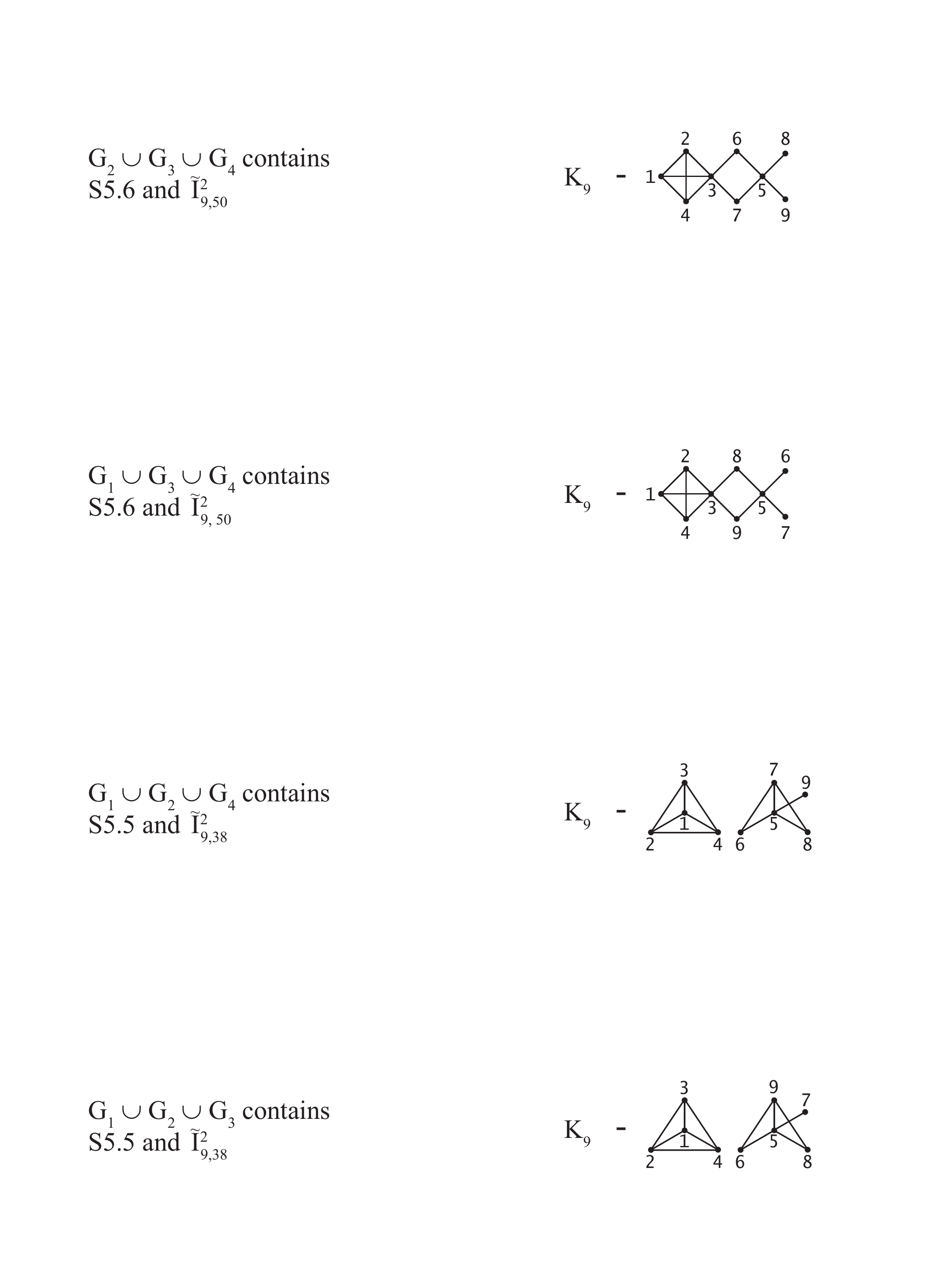}
\caption{\label{fig:kc-v9-n3-9-2} Kuratowski covering of $G = \tilde{I}^3_{9,9}$
for $\mathbb{N}_3$ (Continued)}
\end{figure}

\begin{figure}
\centering
\includegraphics[viewport=0 0 9in 14in, width=12cm]{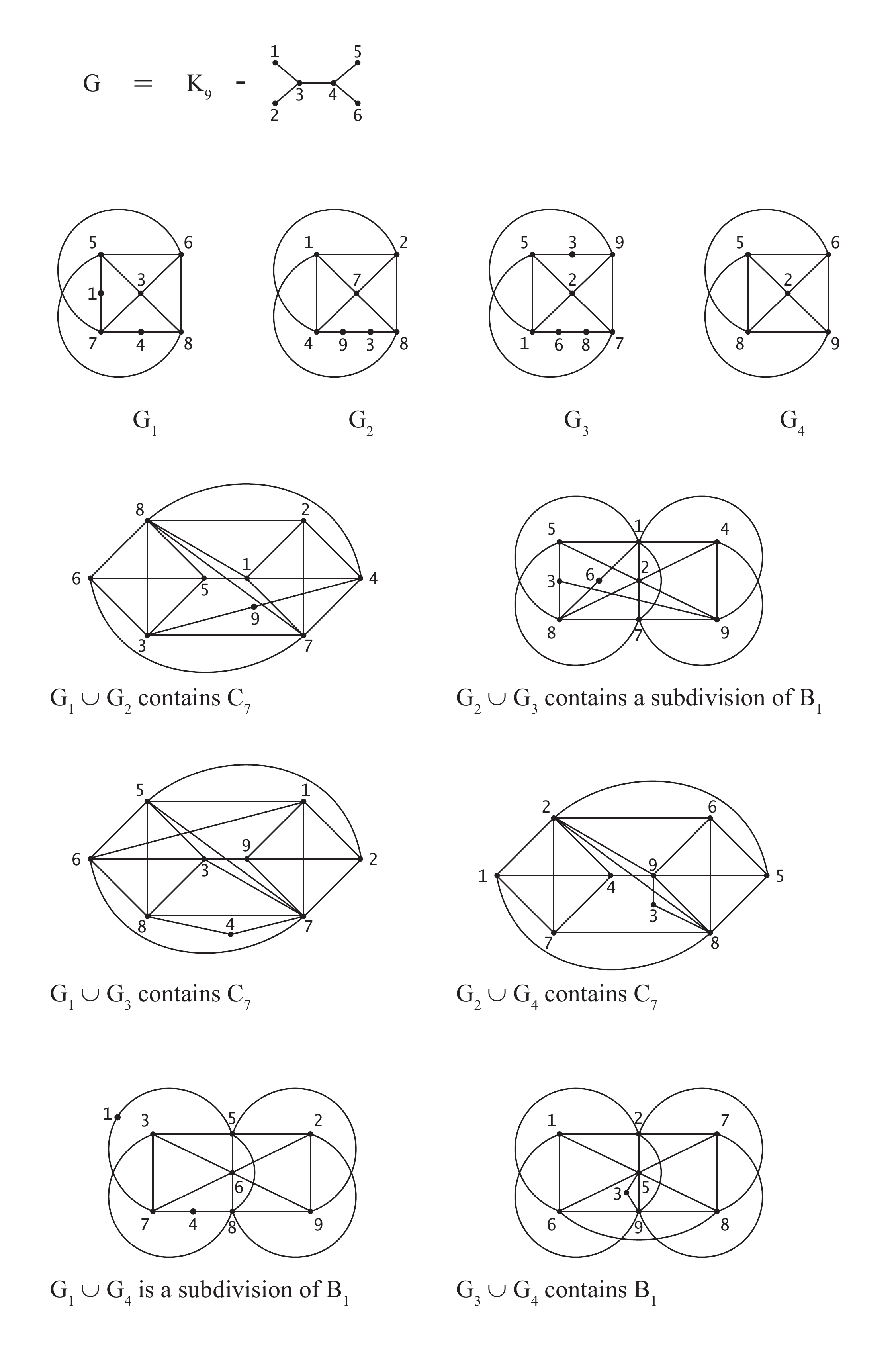}
\caption{\label{fig:kc-v9-n3-10-1} Kuratowski covering of $G = \tilde{I}^3_{9,10}$
for $\mathbb{N}_3$}
\end{figure}

\begin{figure}
\centering
\includegraphics[viewport=0 0 9in 13in, width=12cm]{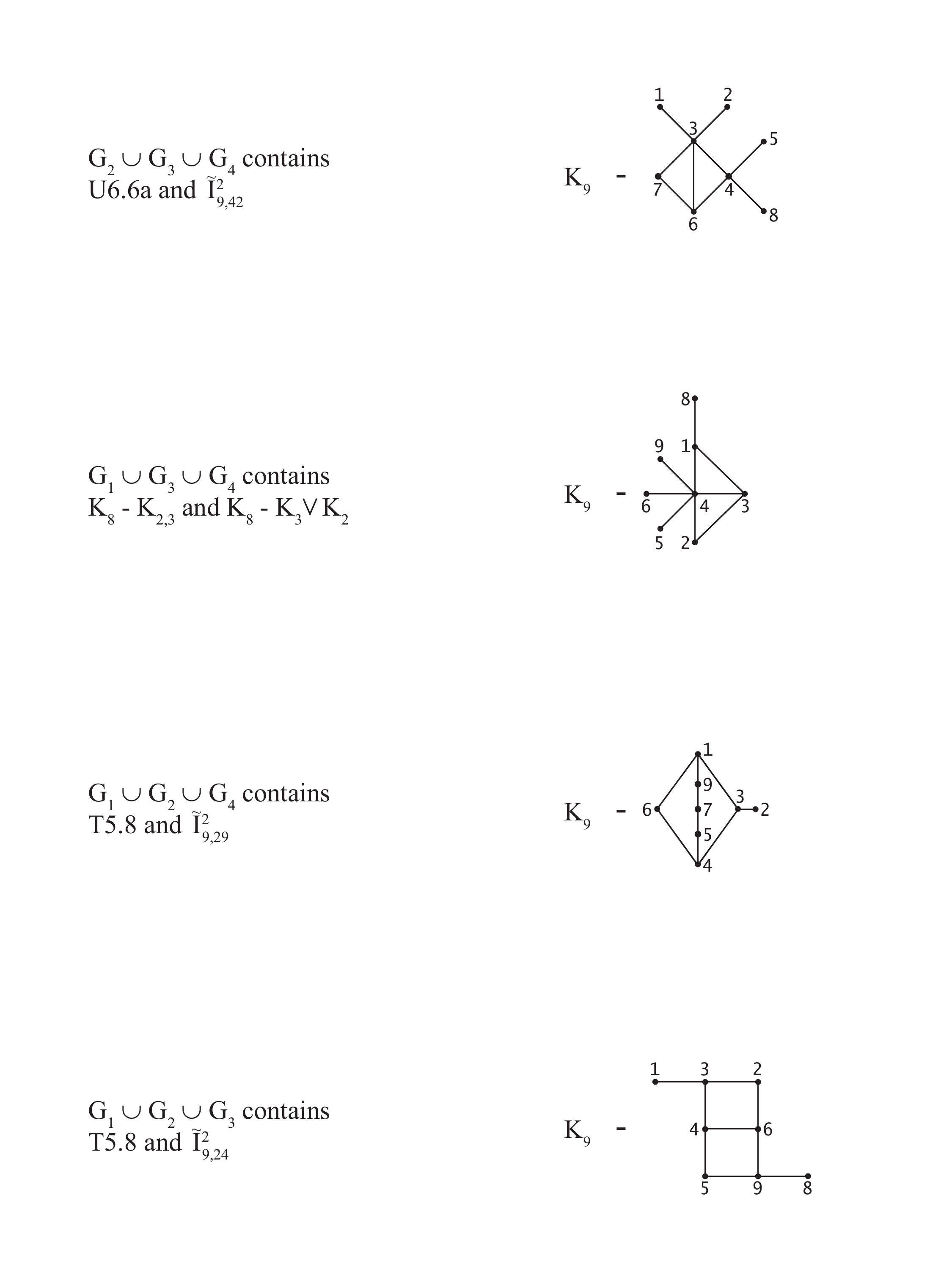}
\caption{\label{fig:kc-v9-n3-10-2} Kuratowski covering of $G = \tilde{I}^3_{9,10}$
for $\mathbb{N}_3$ (Continued)}
\end{figure}

\begin{figure}
\centering
\includegraphics[viewport=0 0 9in 14in, width=12cm]{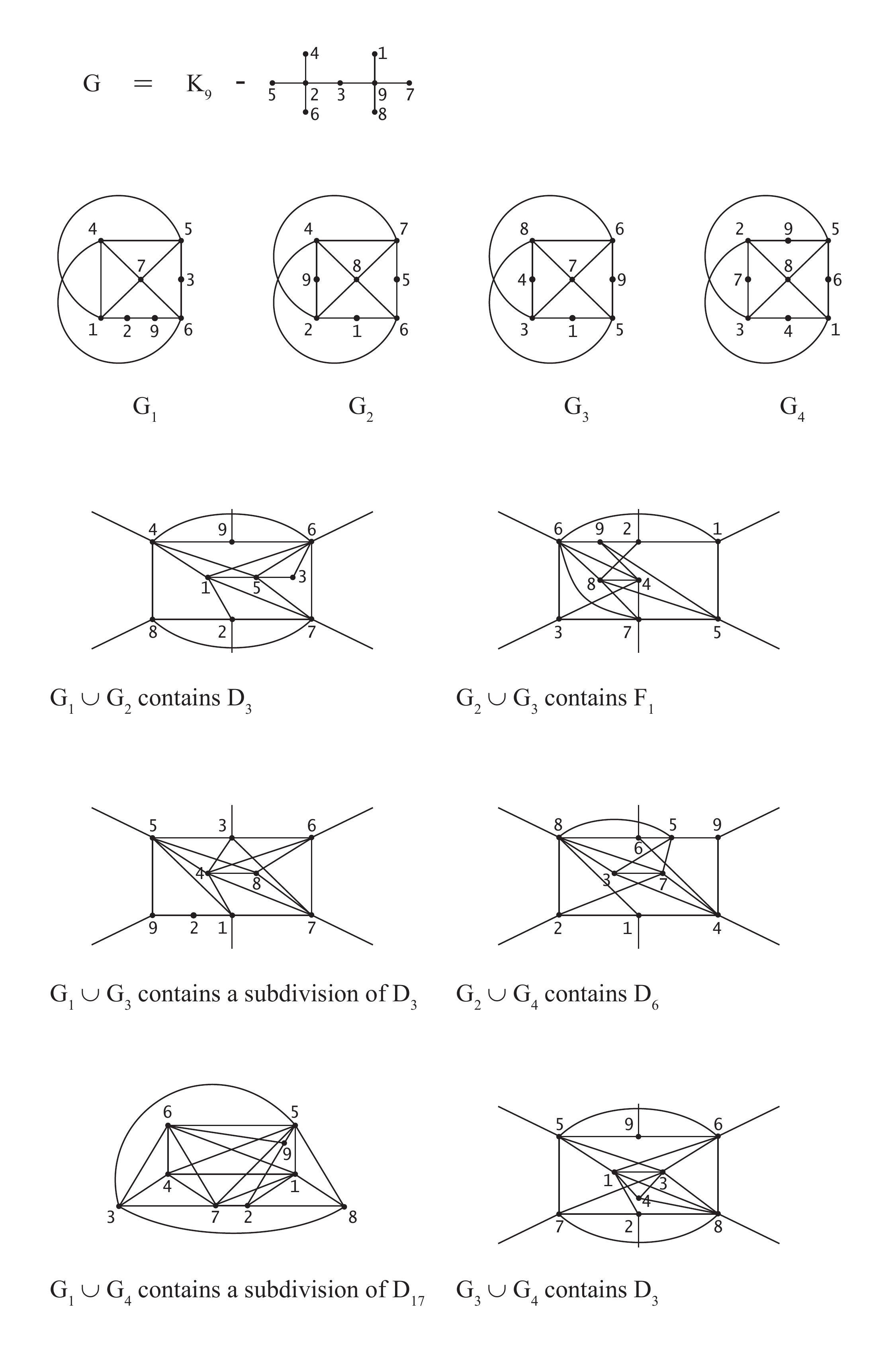}
\caption{\label{fig:kc-v9-n3-11-1} Kuratowski covering of $G = \tilde{I}^3_{9,11}$
for $\mathbb{N}_3$}
\end{figure}

\begin{figure}
\centering
\includegraphics[viewport=0 0 9in 13in, width=12cm]{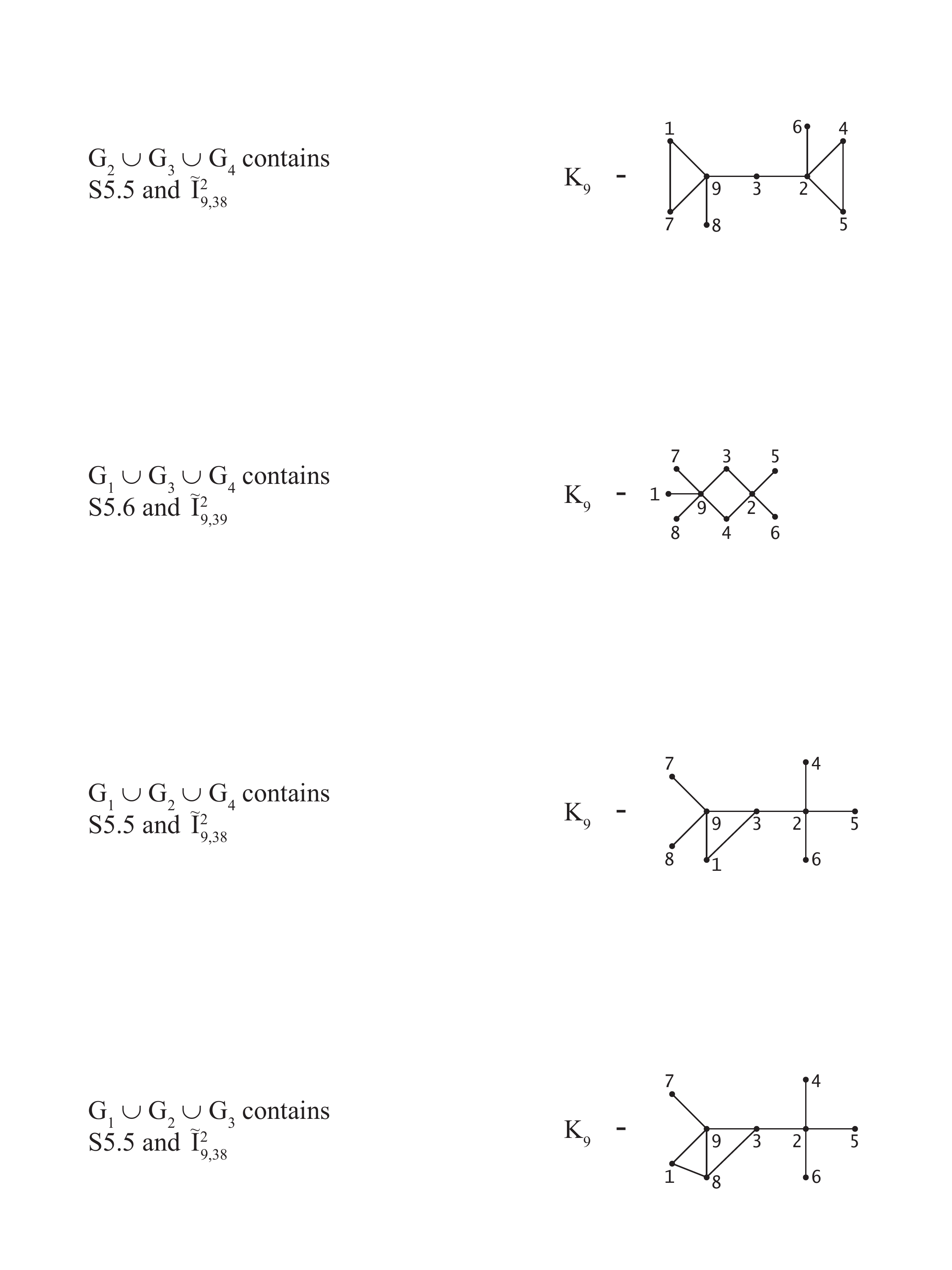}
\caption{\label{fig:kc-v9-n3-11-2} Kuratowski covering of $G = \tilde{I}^3_{9,11}$
for $\mathbb{N}_3$ (Continued)}
\end{figure}

\begin{figure}
\centering
\includegraphics[viewport=0 0 9in 14in, width=12cm]{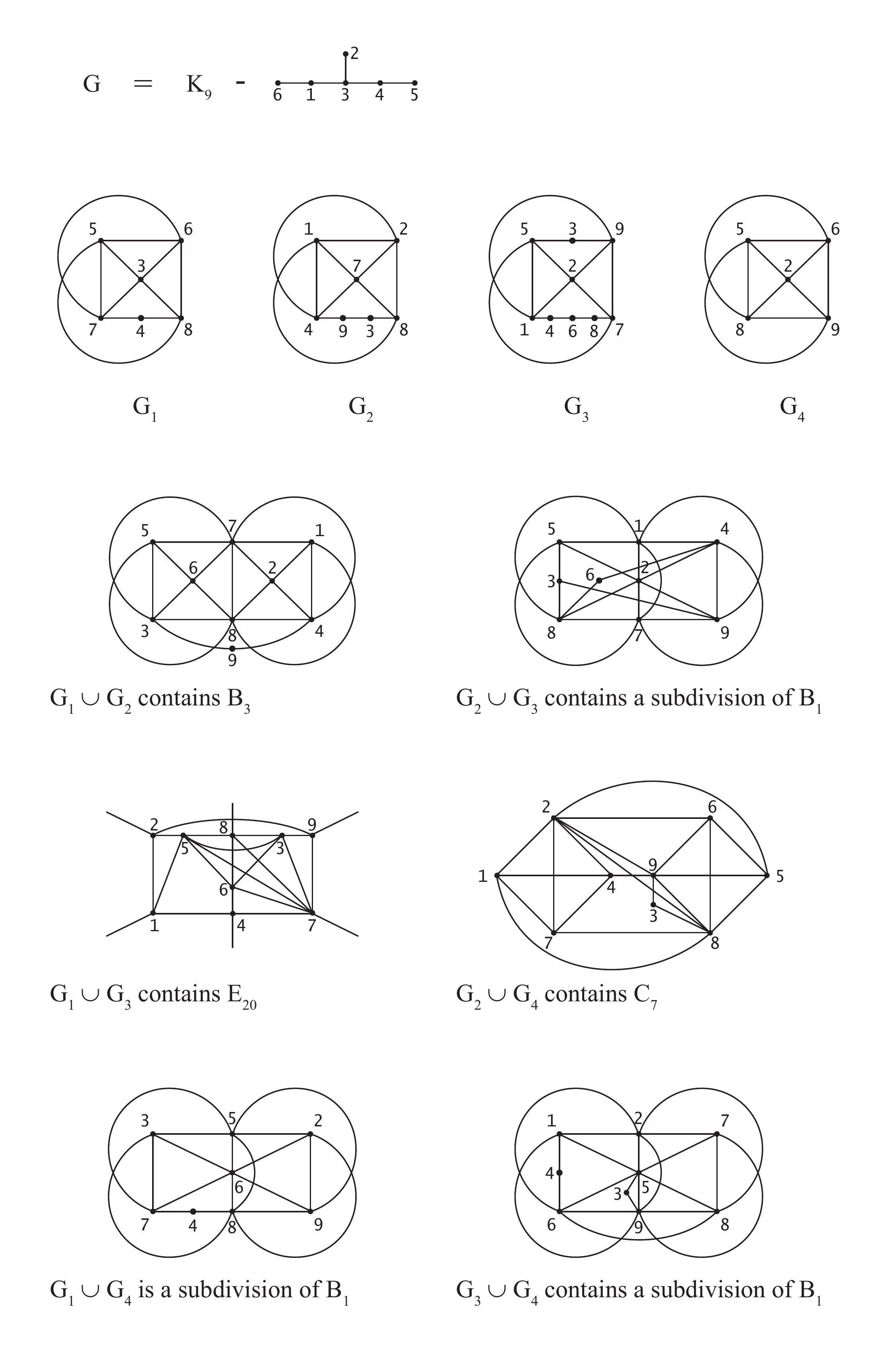}
\caption{\label{fig:kc-v9-n3-12-1} Kuratowski covering of $G = \tilde{I}^3_{9,12}$
for $\mathbb{N}_3$}
\end{figure}

\begin{figure}
\centering
\includegraphics[viewport=0 0 9in 13in, width=12cm]{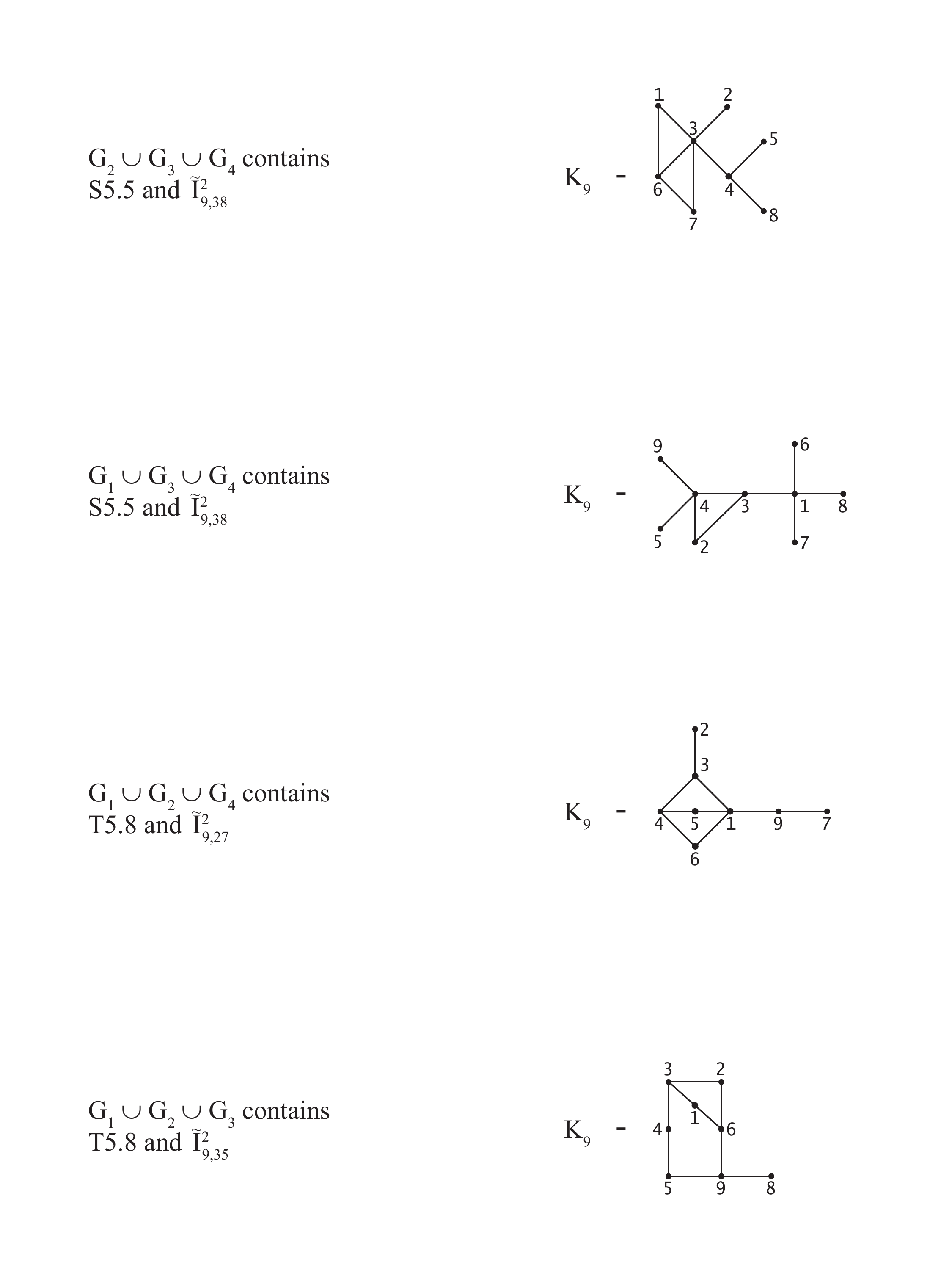}
\caption{\label{fig:kc-v9-n3-12-2} Kuratowski covering of $G = \tilde{I}^3_{9,12}$
for $\mathbb{N}_3$ (Continued)}
\end{figure}

\begin{figure}
\centering
\includegraphics[viewport=0 0 9in 14in, width=12cm]{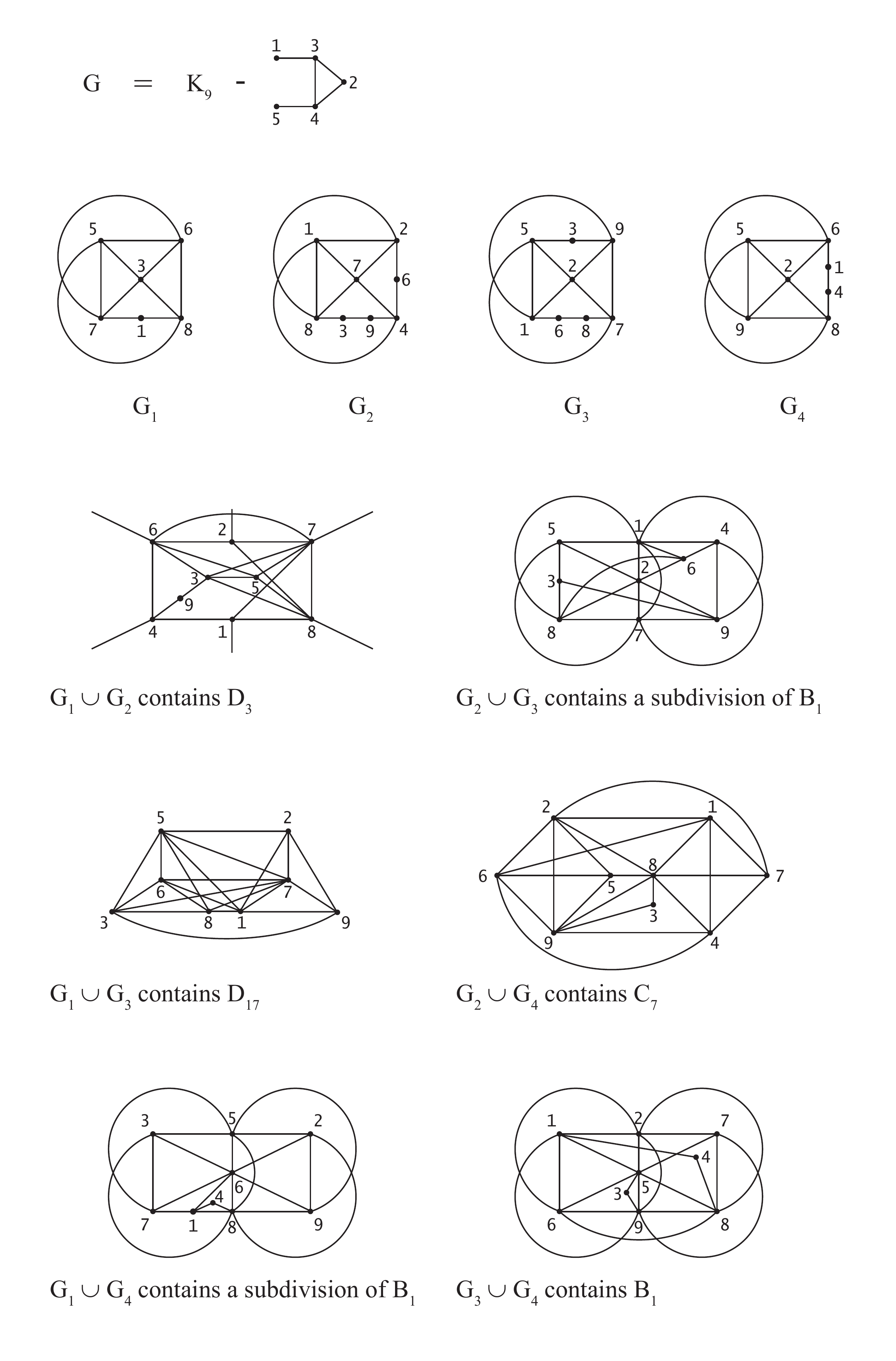}
\caption{\label{fig:kc-v9-n3-13-1} Kuratowski covering of $G = \tilde{I}^3_{9,13}$
for $\mathbb{N}_3$}
\end{figure}

\begin{figure}
\centering
\includegraphics[viewport=0 0 9in 13in, width=12cm]{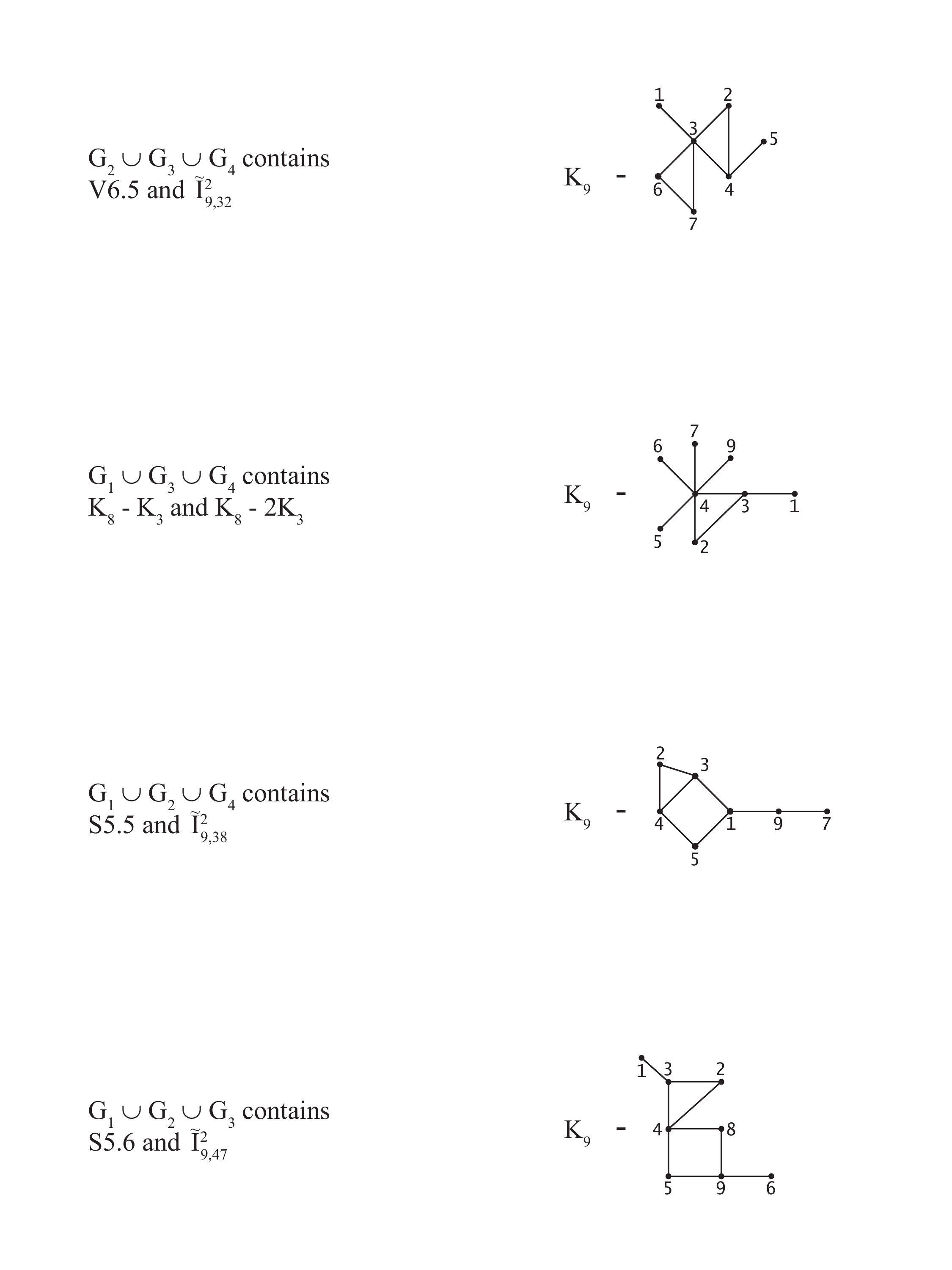}
\caption{\label{fig:kc-v9-n3-13-2} Kuratowski covering of $G = \tilde{I}^3_{9,13}$
for $\mathbb{N}_3$ (Continued)}
\end{figure}

\begin{figure}
\centering
\includegraphics[viewport=0 0 9in 14in, width=12cm]{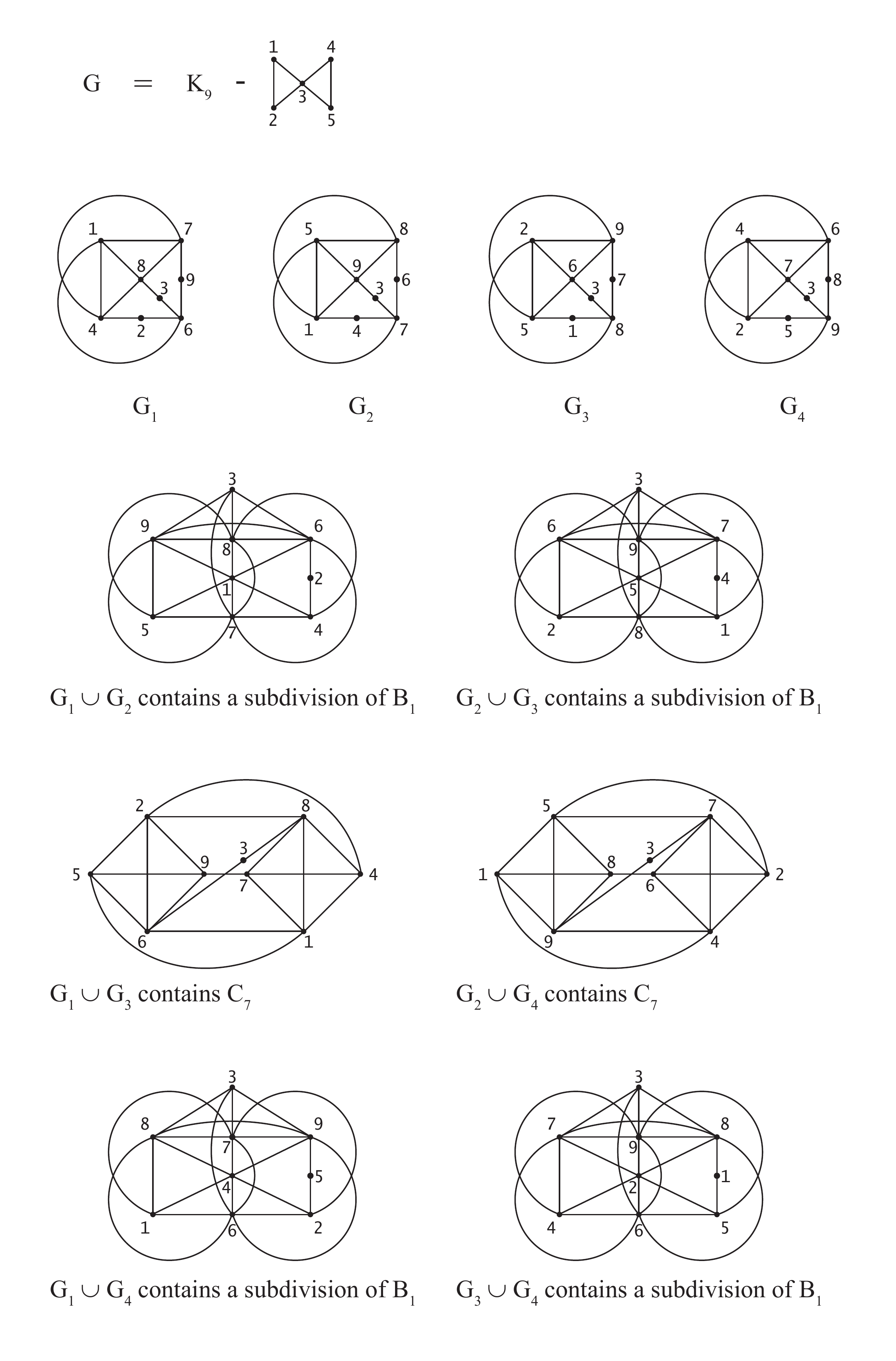}
\caption{\label{fig:kc-v9-n3-14-1} Kuratowski covering of $G = \tilde{I}^3_{9,14}$
for $\mathbb{N}_3$}
\end{figure}

\begin{figure}
\centering
\includegraphics[viewport=0 0 9in 13in, width=12cm]{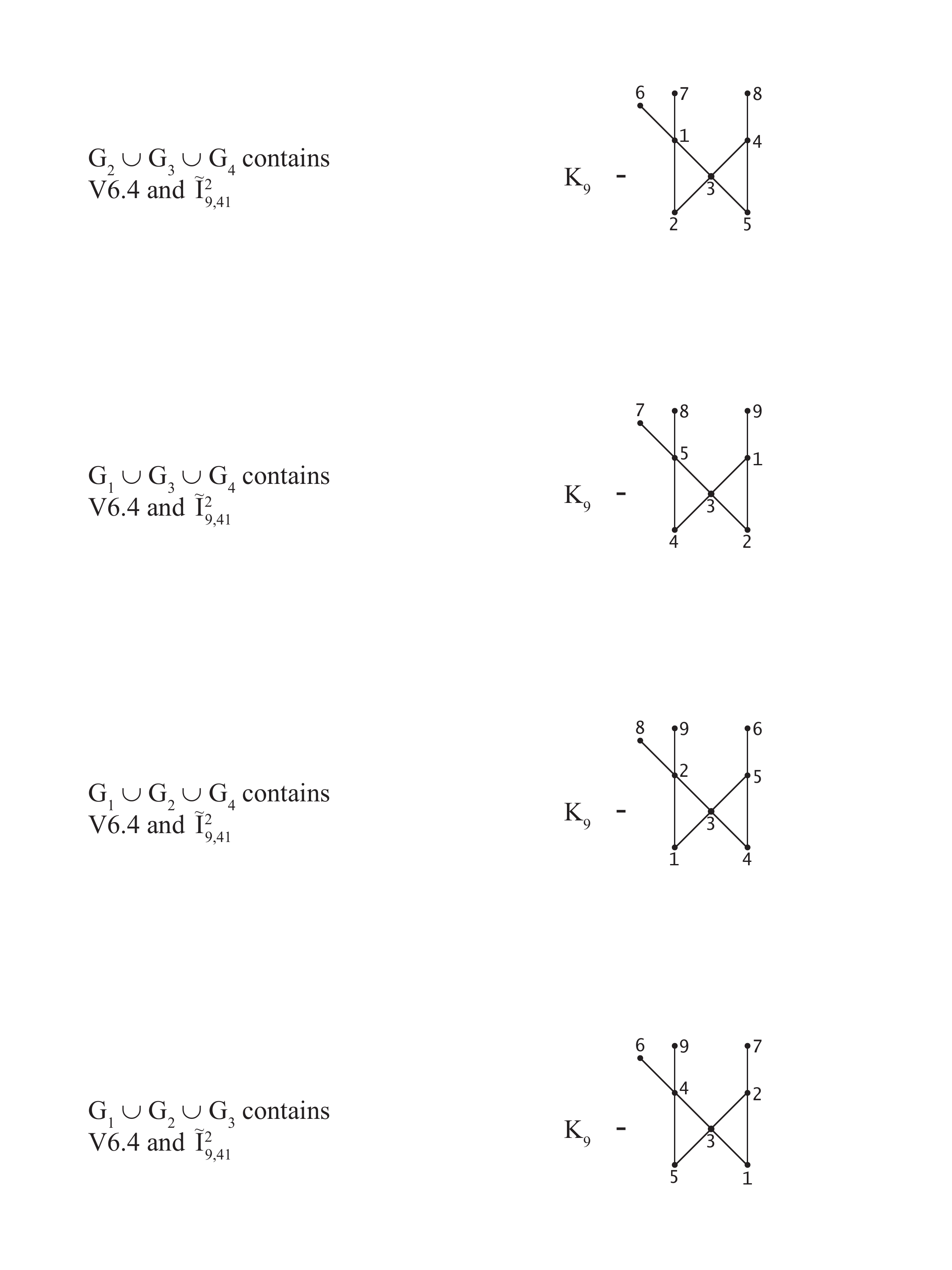}
\caption{\label{fig:kc-v9-n3-14-2} Kuratowski covering of $G = \tilde{I}^3_{9,14}$
for $\mathbb{N}_3$ (Continued)}
\end{figure}

\begin{figure}
\centering
\includegraphics[viewport=0 0 9in 14in, width=12cm]{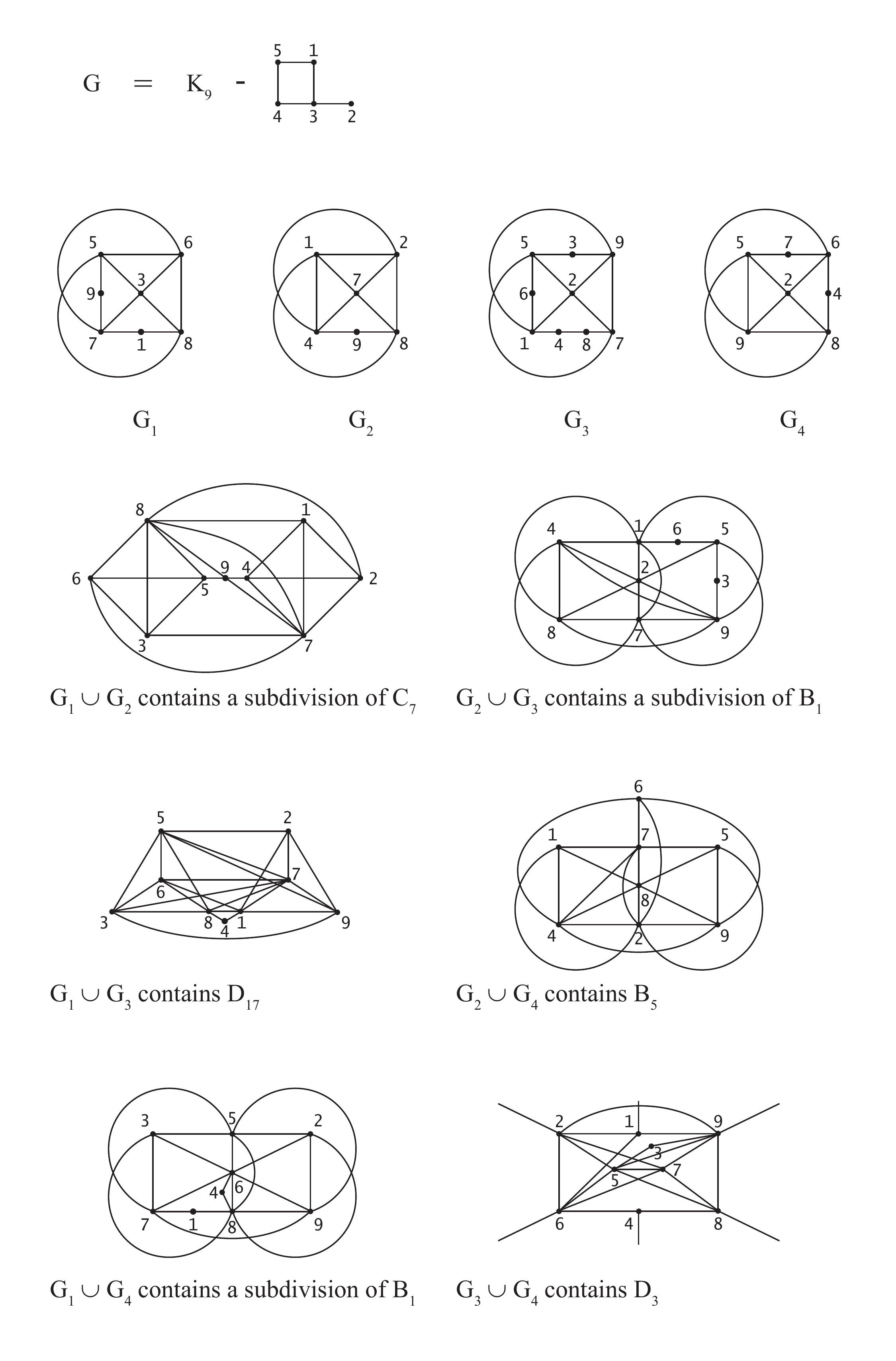}
\caption{\label{fig:kc-v9-n3-15-1} Kuratowski covering of $G = \tilde{I}^3_{9,15}$
for $\mathbb{N}_3$}
\end{figure}

\begin{figure}
\centering
\includegraphics[viewport=0 0 9in 13in, width=12cm]{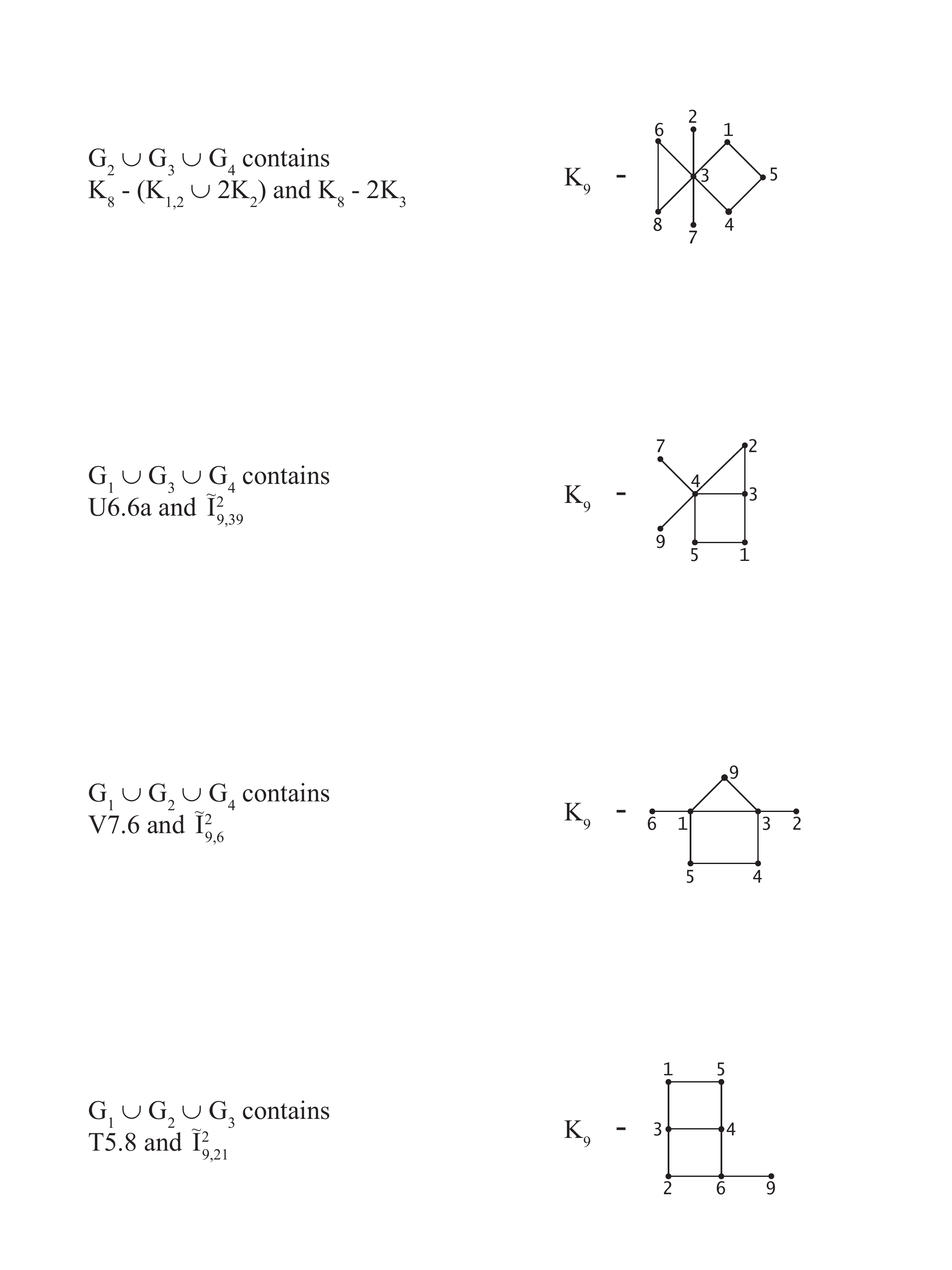}
\caption{\label{fig:kc-v9-n3-15-2} Kuratowski covering of $G = \tilde{I}^3_{9,15}$
for $\mathbb{N}_3$ (Continued)}
\end{figure}

\begin{figure}
\centering
\includegraphics[viewport=0 0 9in 14in, width=12cm]{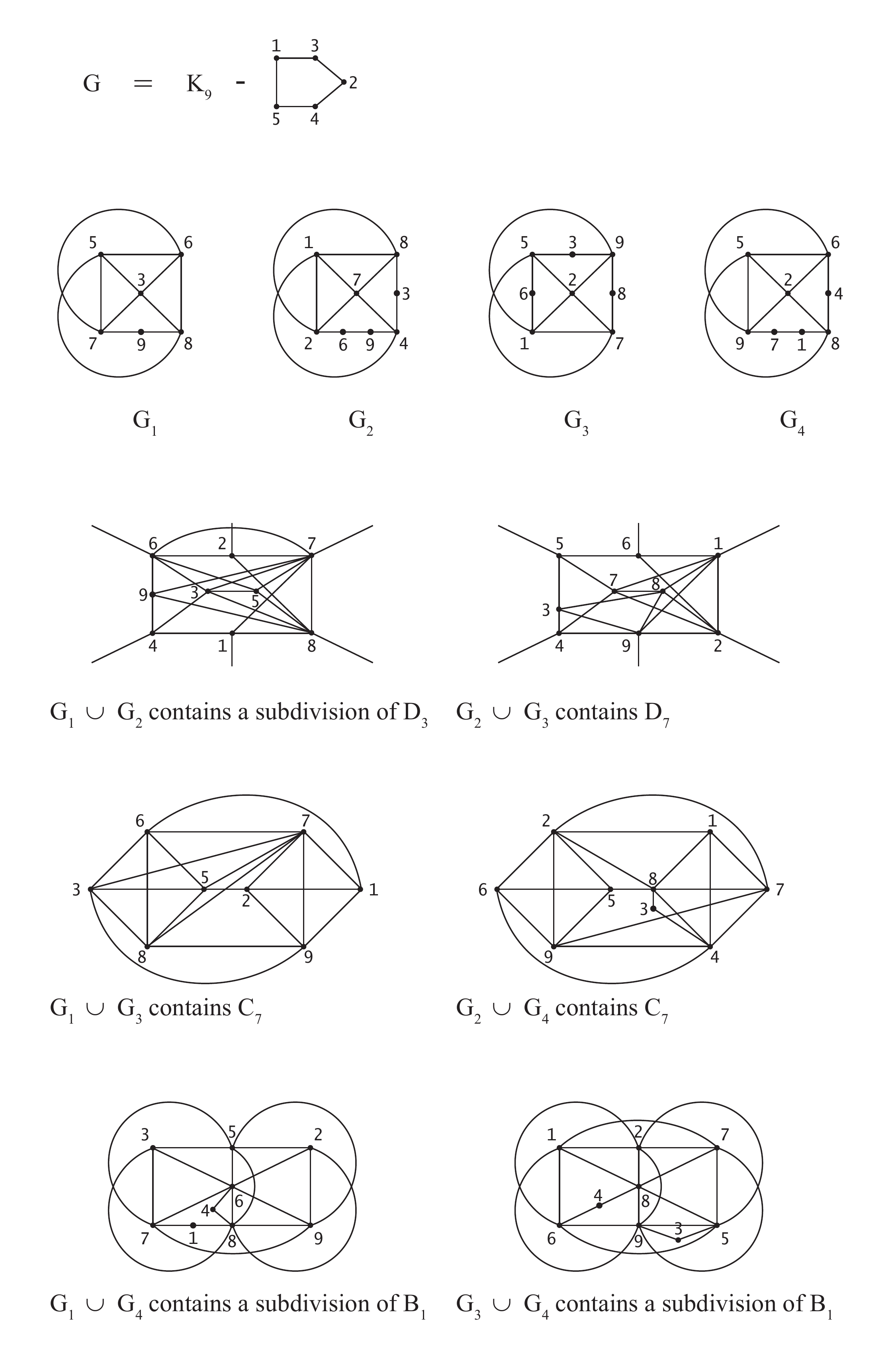}
\caption{\label{fig:kc-v9-n3-16-1} Kuratowski covering of $G = \tilde{I}^3_{9,16}$
for $\mathbb{N}_3$}
\end{figure}

\clearpage

\begin{figure}
\centering
\includegraphics[viewport=0 0 8.5in 14in, width=11cm]{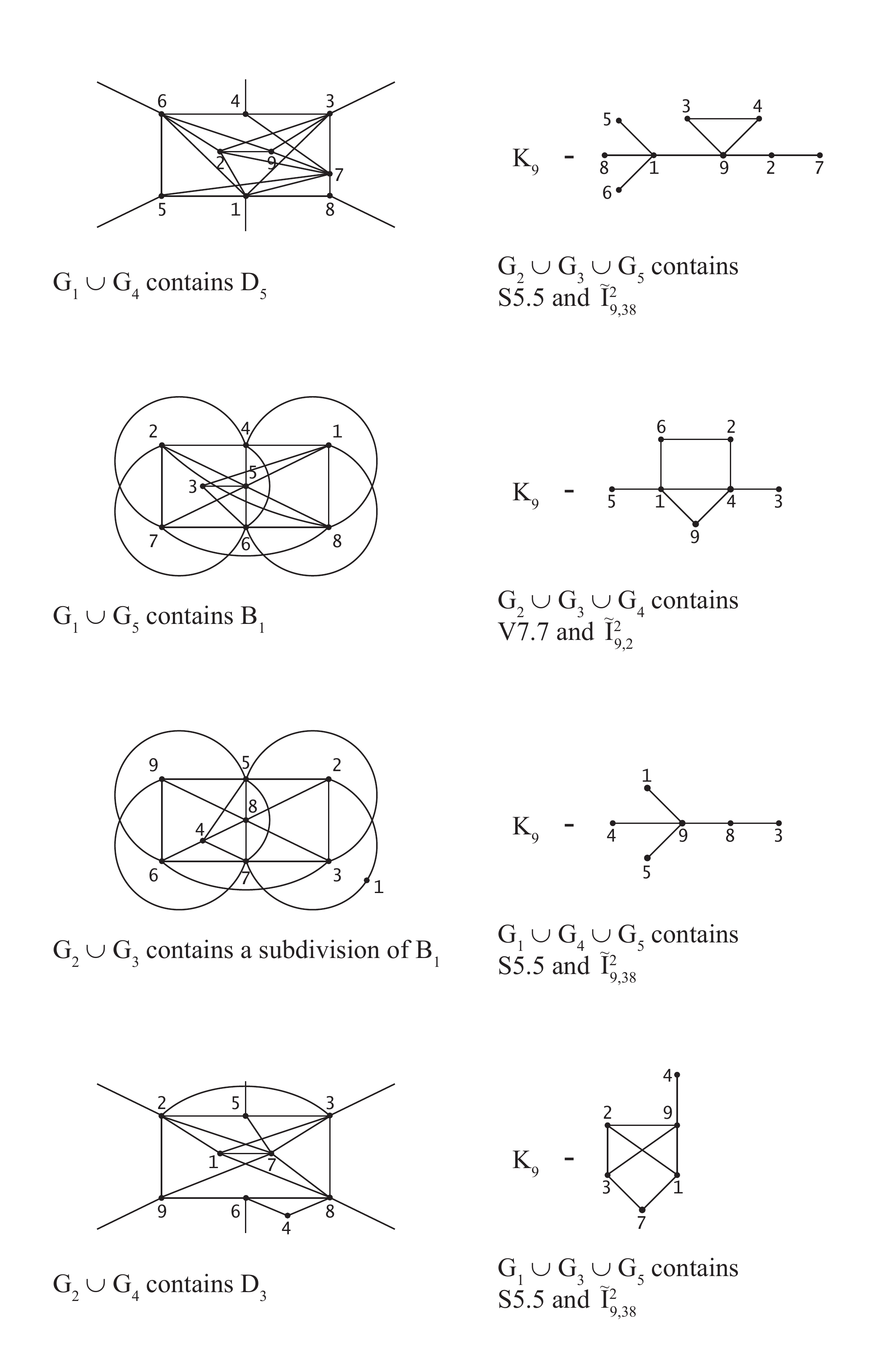}
\caption{\label{fig:kc-v9-n4-1-2} Kuratowski covering of $G = \tilde{I}^4_{9,1}$
for $\mathbb{N}_4$ (Continued)}
\end{figure}

\begin{figure}
\centering
\includegraphics[viewport=0 0 8.5in 14in, width=11cm]{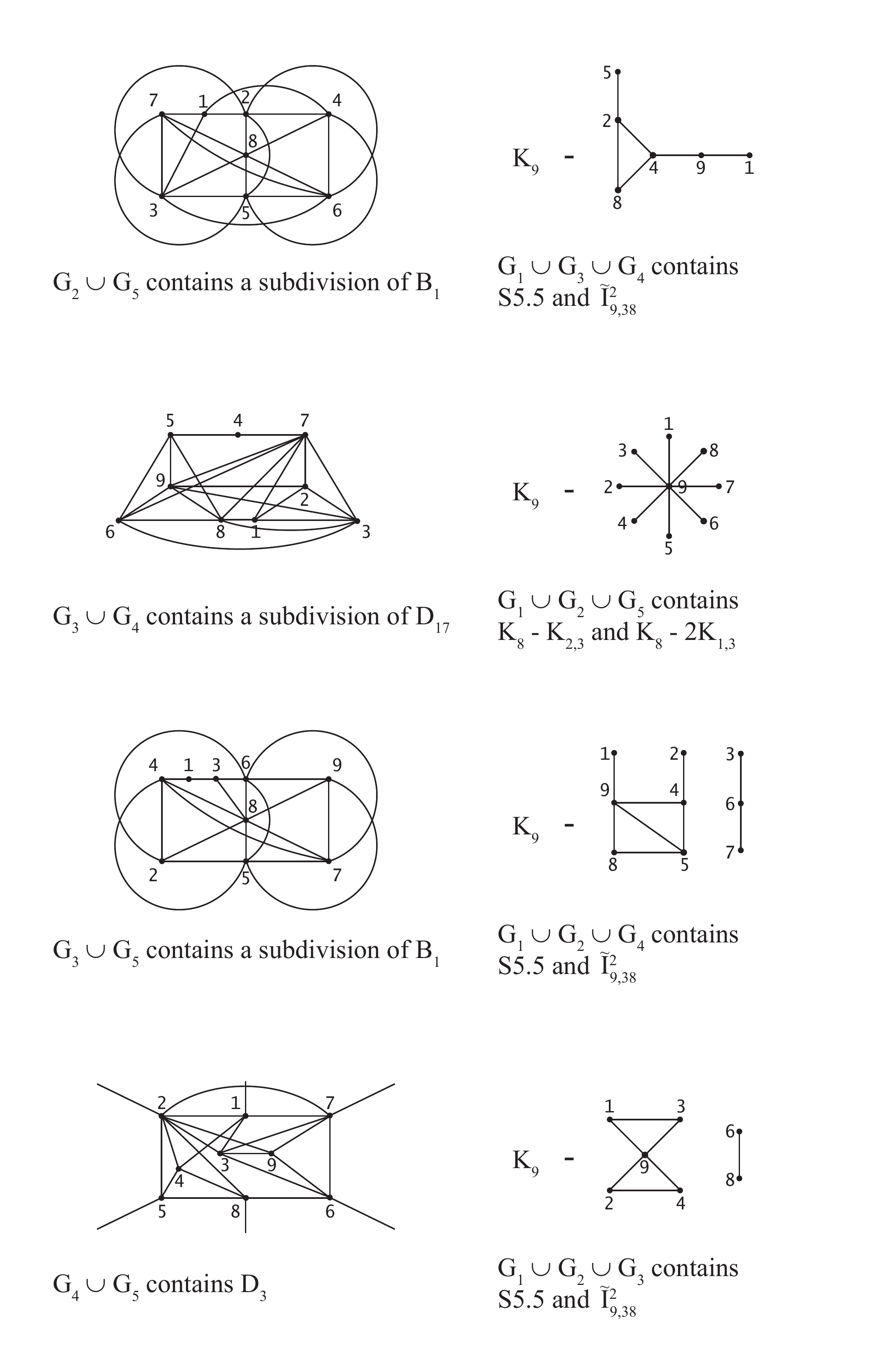}
\caption{\label{fig:kc-v9-n4-1-3} Kuratowski covering of $G = \tilde{I}^4_{9,1}$
for $\mathbb{N}_4$ (Continued)}
\end{figure}

\begin{figure}
\centering
\includegraphics[viewport=0 0 8.5in 14in, width=11cm]{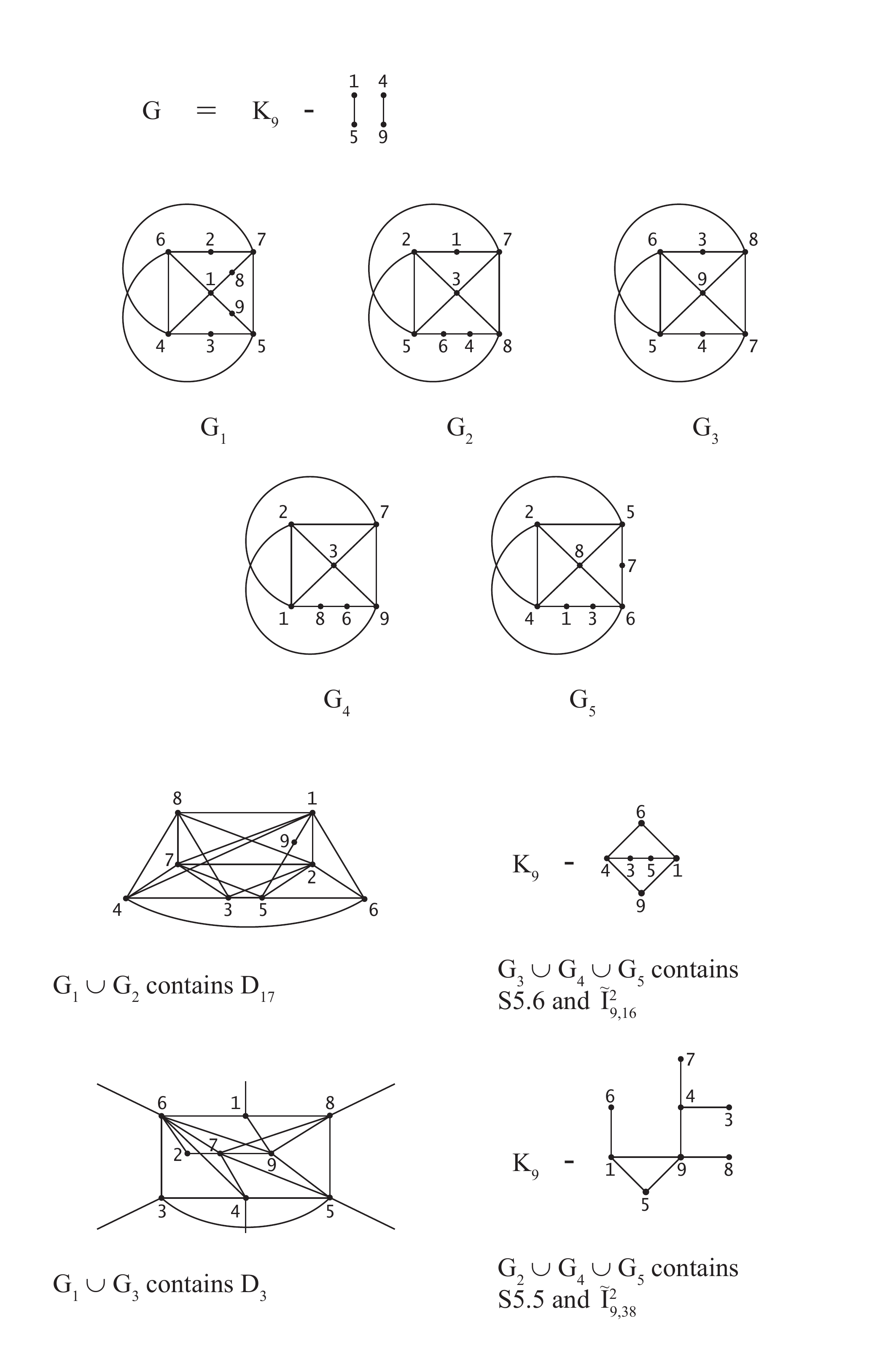}
\caption{\label{fig:kc-v9-n4-2-1} Kuratowski covering of $G = \tilde{I}^4_{9,2}$
for $\mathbb{N}_4$}
\end{figure}

\begin{figure}
\centering
\includegraphics[viewport=0 0 8.5in 14in, width=11cm]{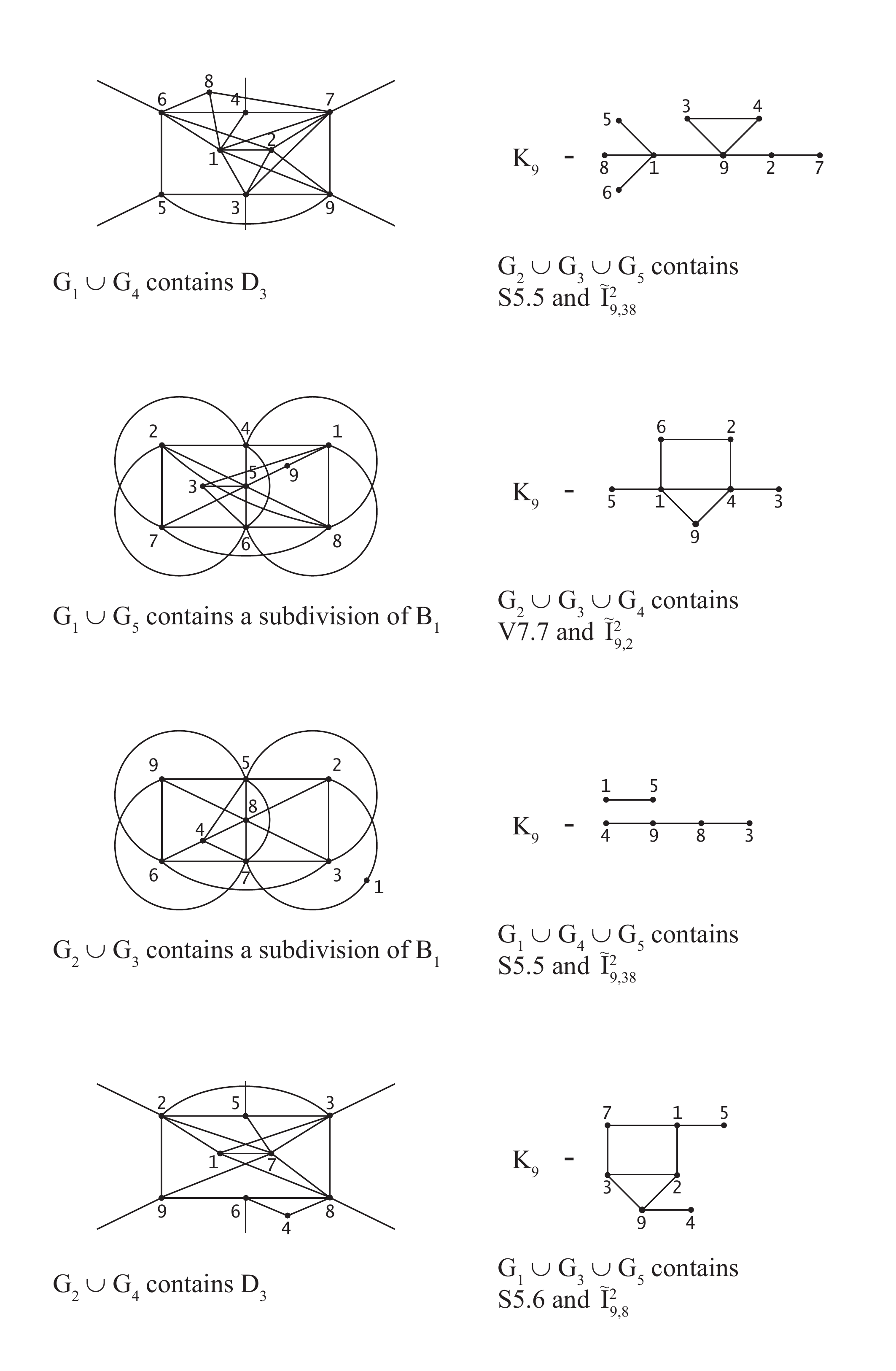}
\caption{\label{fig:kc-v9-n4-2-2} Kuratowski covering of $G = \tilde{I}^4_{9,2}$
for $\mathbb{N}_4$ (Continued)}
\end{figure}

\clearpage

\bibliographystyle{amsplain}
\bibliography{kc-n34}

\end{document}